# TRIMMED TREES AND EMBEDDED PARTICLE SYSTEMS[1]

By Klaus Fleischmann and Jan M. Swart

*Weierstrass Institute and University Erlangen–Nuremberg*

In a supercritical branching particle system, the trimmed tree consists of those particles which have descendants at all times. We develop this concept in the superprocess setting. For a class of continuous superprocesses with Feller underlying motion on compact spaces, we identify the trimmed tree, which turns out to be a binary splitting particle system with a new underlying motion that is a compensated $h$-transform of the old one. We show how trimmed trees may be estimated from above by embedded binary branching particle systems.

**Contents**

Received December 2002; revised July 2003.
[1]Supported by the German Science Foundation, in part by its graduate program "Probabilistic Analysis and Stochastic Processes."

## 1. Introduction and main results.

1.1. *Introduction.* It frequently happens that a superprocess $\mathcal{X} = (\mathcal{X}_t)_{t \geq 0}$, taking values in the space $\mathcal{M}(E)$ of finite measures on some space $E$, and a branching particle process $X = (X_t)_{t \geq 0}$ are related by the formula

$$(1.1) \qquad P^{\mathcal{L}(\text{Pois}(\mu))}[X_t \in \cdot] = P^\mu[\text{Pois}(\mathcal{X}_t) \in \cdot], \qquad t \geq 0, \mu \in \mathcal{M}(E).$$

Here $\text{Pois}(\mathcal{X}_t)$ denotes a Poisson point measure with random intensity $\mathcal{X}_t$ and $P^{\mathcal{L}(\text{Pois}(\mu))}$ denotes the law of the process $X$, started with initial law $\mathcal{L}(X_0) = \mathcal{L}(\text{Pois}(\mu))$. For example, (1.1) holds when $\mathcal{X}$ is the standard, critical, continuous super-Brownian motion in $\mathbb{R}^d$, which corresponds to the evolution equation $\frac{\partial}{\partial t} u_t = \frac{1}{2} \Delta u_t - u_t^2$, and $X$ is a system of binary branching Brownian motions with branching rate 1 and death rate 1. Loosely speaking, $X$ can be obtained from $\mathcal{X}$ by *Poissonization*. Poissonization relations of the form (1.1) have been exploited by various authors, for example, Gorostiza, Roelly-Coppoletta and Wakolbinger ([17], formula (8)), Klenke ([19], formula (4.19)) and Winter ([26], formula (1.23)).

In the present paper, we investigate Poissonization relations for a class of continuous superprocesses on compacta with Feller underlying motion. We give conditions that imply that a superprocess $\mathcal{X}$ and a branching particle system $X$ can be coupled *as processes*, such that

$$(1.2) \qquad P[X_t \in \cdot | (\mathcal{X}_s)_{0 \leq s \leq t}] = P[\text{Pois}(h\mathcal{X}_t) \in \cdot | \mathcal{X}_t] \qquad \text{a.s. } \forall\, t \geq 0,$$

where $h$ is a sufficiently smooth density. Formula (1.2) says that the conditional law of $X_t$, given $(\mathcal{X}_s)_{0 \leq s \leq t}$, is the law of a Poisson point measure with intensity $h\mathcal{X}_t$. For certain critical and subcritical superprocesses, a coupling of the form (1.2) has occurred before in [20], Theorem 3.1 and Section 3.2.

The *weighted superprocess* $(h\mathcal{X}_t)_{t \geq 0}$ that occurs in (1.2) is a superprocess itself, which compared to $\mathcal{X}$ has a new branching mechanism and a new underlying motion, the latter being a "compensated" $h$-transform of the old one. For the special case that $\mathcal{X}$ is a superdiffusion, this fact was proved and exploited by Engländer and Pinsky [7].

Let $\mathcal{X}$ and $X$ be related by (1.2), let $\mathcal{A} := \{\exists \tau < \infty$ such that $\mathcal{X}_t = 0$ $\forall t \geq \tau\}$ denote the event that $\mathcal{X}$ becomes extinct after some random time $\tau$ and set $A := \{\exists \tau < \infty$ such that $X_t = 0\ \forall t \geq \tau\}$. Since $P[X_t = 0 | \mathcal{X}_t = 0] = 1$, $t \geq 0$, we clearly have $P(\mathcal{A} \setminus A) = 0$. We investigate when $\mathcal{X}$ and $X$ can be coupled such that $P(A \setminus \mathcal{A}) = 0$ also holds, that is, the extinction of $X$ implies the extinction of $\mathcal{X}$.

In particular, for a supercritical superprocess $\mathcal{X}$, we construct a binary splitting particle system $X$, such that $\mathcal{X}$ and $X$ are related by a formula of the form (1.2), and, moreover, $X$ corresponds, loosely speaking, to those infinitesimal bits of mass of $\mathcal{X}$ which have descendants at all times. More precisely, we couple the historical processes $\hat{\mathcal{X}}$ and $\hat{X}$ associated with $\mathcal{X}$ and $X$ such that

(1.3) $\forall t \geq 0 \, \exists \tau < \infty$ s.t. $\forall r \geq \tau \qquad \mathrm{supp}(\hat{X}_t) = \mathrm{supp}\,(\hat{\mathcal{X}}_r \circ \pi_{[0,t]}^{-1})$ \qquad a.s.

Here $\pi_{[0,t]}$ denotes projection on the space $\mathcal{D}_E[0,t]$ of cadlag paths from $[0,t]$ into $E$. Informally, $\hat{\mathcal{X}}_t$ is a random measure on paths of length $t$, measuring how much each *line of descent* contributes to the population at time $t$; likewise, $\hat{X}_t$ counts how often each line of descent contributes to $X_t$. Thus, (1.3) says that eventually all mass of the superprocess $\mathcal{X}$ descends from finitely many lines of descent, which are given by $\mathrm{supp}(\hat{X}_t)$. In this special case, the function $h$ that occurrs in (1.2) is $h = p$, the *infinitesimal survival probability* of $\mathcal{X}$, given by

(1.4) $$p(x) = \frac{\partial}{\partial \varepsilon} P^{\varepsilon \delta_x}[\mathcal{X}_t > 0 \; \forall t \geq 0]\bigg|_{\varepsilon=0}, \qquad x \in E.$$

We call $X$ the *trimmed tree* of $\mathcal{X}$. The *reduced tree* of a branching process describes the family relations between all particles alive at a fixed time and their ancestors (neglecting those lines of descent that died earlier). Thus, our trimmed tree can be viewed as the limit of reduced trees as time tends to infinity. Reduced trees have been studied intensively in the branching literature. For historical background, see, for example, the last paragraph in Section 12.1 of [2], page 201.

It is worth mentioning that the weighted superprocess $(p\mathcal{X}_t)_{t \geq 0}$ with $p$ as in (1.4) played an important role in the work of Engländer and Pinsky [7], who investigated support properties (such as recurrence) of superdiffusions by analytic tools. Weighted superprocesses and embedded particle systems also played a central role in [14], which motivated our present article.

The paper is organized as follows. In Sections 1.2–1.4, we introduce our objects of interest together with some of their elementary properties in more detail. Sections 1.5 and 1.6 contain our main results, while Section 1.7 is devoted to discussion. In Section 2, we collect some necessary facts on historical processes and weighted superprocesses. The final proofs are deferred to Section 3.

1.2. *Poissonization of superprocesses.* Let $E$ be a compact metrizable space, and let $B(E)$ and $\mathcal{C}(E)$ denote the spaces of bounded measurable real functions and continuous real functions on $E$, respectively. We set $B_+(E) := \{f \in B(E) : f \geq 0\}$, $B_{[0,1]}(E) := \{f \in B(E) : 0 \leq f \leq 1\}$ and define $\mathcal{C}_+(E)$, $\mathcal{C}_{[0,1]}(E)$ similarly. Let $\mathcal{M}(E)$ denote the space of finite measures on $E$, equipped with the topology of weak convergence. If $\mu$ is a measure and $f$ is measurable, then $\langle \mu, f \rangle := \int_E f \, d\mu$ denotes the integral of $f$ with respect to $\mu$, whenever it exists. By $\mathcal{N}(E) \subset \mathcal{M}(E)$ we denote the space of finite point measures, that is, measures $\nu$ of the form $\sum_{i=1}^{n} \delta_{x_i}$ with $x_i \in E$ and $n \geq 0$. We interpret such a point measure as a collection of $n$ particles, situated at positions $x_1, \ldots, x_n$. For $f \in B_{[0,1]}(E)$ and $\nu = \sum_{i=1}^{n} \delta_{x_i} \in \mathcal{N}(E)$ we use the notation $f^\nu := \prod_{i=1}^{n} f(x_i)$ (where $f^0 := 1$). If $\mu$ is a random variable taking values in $\mathcal{M}(E)$, then $\text{Pois}(\mu)$ denotes an $\mathcal{N}(E)$-valued random variable such that, conditioned on $\mu$, $\text{Pois}(\mu)$ is a Poisson point measure with intensity $\mu$. If $\nu$ is a random variable taking values in $\mathcal{N}(E)$ and $f \in B_{[0,1]}(E)$, then $\text{Thin}_f(\nu)$ denotes a random point measure obtained by thinning $\nu$ with $f$. That is, conditioned on $\nu$, a particle $\delta_x$ in $\nu$ is kept with probability $f(x)$, independently of the other particles in $\nu$. Note that

(1.5)
$$\begin{aligned}&\text{(i) } P[\text{Pois}(f\mu) = 0|\mu] = e^{-\langle \mu, f \rangle}, \qquad f \in B_+(E),\\ &\text{(ii) } P[\text{Thin}_f(\nu) = 0|\nu] = (1-f)^\nu, \quad f \in B_{[0,1]}(E).\end{aligned}$$

It is well known that

(1.6) $\text{Thin}_f(\text{Thin}_g(\nu)) \stackrel{\mathcal{D}}{=} \text{Thin}_{fg}(\nu)$ and $\text{Thin}_f(\text{Pois}(\mu)) \stackrel{\mathcal{D}}{=} \text{Pois}(f\mu)$,

where $\stackrel{\mathcal{D}}{=}$ denotes equality in distribution.

Let $G$ be the generator of a Feller process $\xi = (\xi_t)_{t \geq 0}$ on $E$ and let $\alpha \in \mathcal{C}_+(E)$, $\beta \in \mathcal{C}(E)$. Then, for each $f \in B_+(E)$, an appropriate integrated version [see (2.8)] of the semilinear Cauchy problem

(1.7)
$$\frac{\partial}{\partial t} u_t = G u_t + \beta u_t - \alpha u_t^2, \qquad t \geq 0,$$
$$u_0 = f,$$

has a unique solution $u_t =: \mathcal{U}_t f$, $t \geq 0$, in $B_+(E)$. Moreover, there exists a unique (in law) Markov process $\mathcal{X}$ in $\mathcal{M}(E)$ with continuous sample paths, defined by its Laplace functionals

(1.8) $\qquad E^\mu[e^{-\langle \mathcal{X}_t, f \rangle}] = e^{-\langle \mu, \mathcal{U}_t f \rangle}, \qquad t \geq 0, \mu \in \mathcal{M}(E), f \in B_+(E).$

The process $\mathcal{X}$ is called the *superprocess* in $E$ with *underlying motion generator* $G$, (local) *activity* $\alpha$ and (local) *growth parameter* $\beta$ (the last two terms are our terminology) or, for short, the $(G, \alpha, \beta)$-*superprocess*. The semigroup $(\mathcal{U}_t)_{t \geq 0} = \mathcal{U} = \mathcal{U}(G, \alpha, \beta)$ is called the *log-Laplace semigroup* of $\mathcal{X}$. In fact, $\mathcal{U}_t f$ can be defined unambiguously for any measurable $f : E \to [0, \infty]$ such that (1.8) holds (where $e^{-\infty} := 0$). The process $\mathcal{X}$ can be constructed in several ways and is nowadays standard (see, e.g., [[10]–[12]]). We can think of $\mathcal{X}$ as describing a population where mass flows with generator $G$ and during a time interval $dt$ a bit of mass $dm$ at position $x$ produces offspring with mean $(1 + \beta(x) dt) dm$ and finite variance $2\alpha(x) dt dm$. For basic facts on superprocesses, we refer to [2] and [8].

Similarly, when $G$ is (again) the generator of a Feller process in a compact metrizable space $E$ and $b, d \in \mathcal{C}_+(E)$, then, for any $f \in B_{[0,1]}(E)$, an integrated version of the semilinear Cauchy problem

(1.9)
$$\frac{\partial}{\partial t} u_t = G u_t + b u_t (1 - u_t) - d u_t, \qquad t \geq 0,$$
$$u_0 = f,$$

has a unique solution $u_t =: U_t f$, $t \geq 0$, in $B_{[0,1]}(E)$. Moreover, there exists a unique Markov process $X$ with cadlag sample paths in $\mathcal{N}(E)$, defined by its generating functionals

(1.10) $\quad E^\nu[(1-f)^{X_t}] = (1 - U_t f)^\nu, \qquad t \geq 0, \nu \in \mathcal{N}(E), f \in B_{[0,1]}(E).$

We call $X$ the *binary branching particle system* in $E$ with underlying motion generator $G$, branching rate $b$ and death rate $d$, or, for short the $(G, b, d)$-*particle system*. The semigroup $(U_t)_{t \geq 0} = U = U(G, b, d)$ is called the *generating semigroup* of $X$. The particles in $X$ perform independent motions with generator $G$ and, additionally, a particle branches with local rate $b$ into two new particles, created at the position of the old one, and particles die with local rate $d$. If the death rate is zero, we also speak about *binary splitting* instead of binary branching.

Because of (1.5), formulas (1.8) and (1.10) can be rewritten as

(1.11) 
(i) $P^\mu[\text{Pois}(f \mathcal{X}_t) = 0] = P[\text{Pois}((\mathcal{U}_t f)\mu) = 0], \mu \in \mathcal{M}(E), f \in B_+(E),$
(ii) $P^\nu[\text{Thin}_f(X_t) = 0] = P[\text{Thin}_{\mathcal{U}_t f}(\nu) = 0], \nu \in \mathcal{N}(E), f \in B_{[0,1]}(E),$

$t \geq 0$. The following lemma is now an easy observation.

LEMMA 1 (Poissonization of superprocesses). *Let $\mathcal{X}$ be the $(G, \alpha, \beta)$-superprocess, assume that $\alpha \geq \beta$ and let $X$ be the $(G, \alpha, \alpha - \beta)$-particle system. Then*

(1.12) $\quad P^{\mathcal{L}(\text{Pois}(\mu))}[X_t \in \cdot] = P^\mu[\text{Pois}(\mathcal{X}_t) \in \cdot], \qquad t \geq 0, \mu \in \mathcal{M}(E).$

PROOF. Let $\mathcal{U} = \mathcal{U}(G, \alpha, \beta)$ and $U = U(G, \alpha, \alpha - \beta)$ denote the log-Laplace semigroup of $\mathcal{X}$ and the moment generating semigroup of $X$, respectively. Comparing the Cauchy problems (1.7) and (1.9), we see that $\mathcal{U}_t f = U_t f$ for all $f \in B_{[0,1]}(E)$ and $t \geq 0$. It follows that for any $f \in B_{[0,1]}(E)$, $\mu \in \mathcal{M}(E)$ and $t \geq 0$,

(1.13)
$$P^{\mathcal{L}(\text{Pois}(\mu))}[\text{Thin}_f(X_t) = 0]$$
$$= P[\text{Thin}_{U_t f}(\text{Pois}(\mu)) = 0] = P[\text{Pois}((U_t f)\mu) = 0]$$
$$= P^\mu[\text{Pois}(f \mathcal{X}_t) = 0] = P^\mu[\text{Thin}_f(\text{Pois}(\mathcal{X}_t)) = 0].$$

Since this holds for arbitrary $f \in B_{[0,1]}(E)$, the law of $X_t$, when $X$ is started with initial law $\mathcal{L}(X_0) = \mathcal{L}(\text{Pois}(\mu))$, coincides with the law of $\text{Pois}(\mathcal{X}_t)$, when $\mathcal{X}$ is started in $\mathcal{X}_0 = \mu$. □

REMARK 2 (Locally compact spaces). With the help of a suitable compactification, the results in this paper can be applied to superprocesses on some noncompact spaces as well. Let $E$ be a locally compact but not compact, separable, metrizable space, let $G$ be the generator of a Feller process $\xi = (\xi_t)_{t \geq 0}$ on $E$, whose semigroup maps the space $\mathcal{C}_0(E)$ of continuous real functions vanishing at infinity into itself and let $\alpha, \beta$ be bounded continuous functions on $E$, $\alpha \geq 0$. Then the $(G, \alpha, \beta)$-superprocess may be defined as follows. First, $E$ may be embedded in a compact metrizable space $\overline{E}$ such that $E$ is an open dense subset of $\overline{E}$ and such that the functions $\alpha, \beta$ can be extended to continuous functions $\overline{\alpha}, \overline{\beta}$ on $\overline{E}$. [To construct such a compactification, take for $\overline{E}$ the closure of the graph of $(\alpha, \beta)$ in $E_\sigma \times \mathbb{R}^2$, where $E_\sigma$ is the one-point compactification of $E$.] Second, $\xi$ may be extended to a Feller process in $\overline{E}$ (with generator denoted by $\overline{G}$) by putting $P^x[\xi_t = x \, \forall t \geq 0] := 1$ for $x \in \overline{E} \setminus E$. By identifying $\mathcal{M}(E)$ with the space $\{\mu \in \mathcal{M}(\overline{E}) : \mu(\overline{E} \setminus E) = 0\}$, the $(\overline{G}, \overline{\alpha}, \overline{\beta})$-superprocess $\overline{\mathcal{X}}$ satisfies $P^\mu[\overline{\mathcal{X}}_t \in \mathcal{M}(E) \; \forall t \geq 0] = 1$ for all $\mu \in \mathcal{M}(E)$. The $(G, \alpha, \beta)$-superprocess may then be defined as the restriction of $\overline{\mathcal{X}}$ to $\mathcal{M}(E)$. In this way, the results in this paper can be applied, for example, to the usual super-Brownian motion (with finite initial mass). To keep notation simple, we formulate our results in the rest of this paper for superprocesses in a compact space $E$.

1.3. *Historical superprocesses and branching particle systems.* Let $E$ be a compact metrizable space as before, and let $\mathcal{D}_E[0, \infty)$ and $\mathcal{D}_E[0, t]$ denote the spaces of cadlag paths $w : [0, \infty) \to E$ and $w : [0, t] \to E$, respectively, equipped with the Skorohod topology. Let $\xi$ be a Feller process in $E$. Then the *path process* $\hat{\xi}$ associated with $\xi$ is a time-inhomogeneous Markov process with time-dependent state space $\mathcal{D}_E[0, t]$, defined as follows. Let $\xi^x$ denote the process $\xi$ started in $\xi_0^x = x \in E$. Then $(\hat{\xi}_t^{s,w})_{t \geq s}$, the path process $\hat{\xi}$ started at time $s \geq 0$ in $w \in \mathcal{D}_E[0, s]$ and evaluated at times $t \geq s$, is defined as

$$(1.14) \qquad \hat{\xi}_t^{s,w}(r) := \begin{cases} w(r), & \text{if } 0 \leq r \leq s, \\ \xi_{r-s}^{w(s)}, & \text{if } s \leq r \leq t. \end{cases}$$

For $t \geq 0$, we identify the space $\mathcal{D}_E[0, t]$ with the space $\{w \in \mathcal{D}_E[0, \infty) : w(u) = w(t) \, \forall u \geq t\}$ of paths stopped at time $t$. With this identification, $\hat{\xi}^{s,w} : [s, \infty) \to \mathcal{D}_E[0, \infty)$ has cadlag sample paths. Note that $\hat{\xi}_t^{0,x}$, the path process started at time zero in $x \in \mathcal{D}_E\{0\} \cong E$ and evaluated at time $t \geq 0$, records the path followed by $\xi^x$ up to time $t$.

If $\mathcal{X}$ is a $(G, \alpha, \beta)$-superprocess in $E$ as defined in the last section, then by definition the *historical superprocess* $\hat{\mathcal{X}}$ associated with $\mathcal{X}$ is the time-inhomogeneous superprocess with time-dependent state space $\mathcal{M}(\mathcal{D}_E[0, t])$, with underlying motion $\hat{\xi}$, time-dependent activity $\hat{\alpha}_t(w) := \alpha(w(t))$ and time-dependent growth parameter $\hat{\beta}_t(w) := \beta(w(t))$. We call $\hat{\mathcal{X}}$ the *historical* $(G, \alpha, \beta)$-*superprocess*. As before, we identify $\mathcal{D}_E[0, t]$ with the subspace of $\mathcal{D}_E[0, \infty)$ consisting of paths stopped at time $t$, and in this identification $\mathcal{X}:[0, \infty) \to \mathcal{M}(\mathcal{D}_E[0, \infty))$ has continuous sample paths. For the technical details needed to deal with the facts that the underlying motion is time-inhomogeneous and the space $\mathcal{D}_E[0, \infty)$ is not locally compact, we refer to Section 2.2; see also [3], Chapter 2. If $\hat{\mathcal{X}}$ is started at time zero in $\hat{\mathcal{X}}_0 = \mu \in \mathcal{M}(\mathcal{D}_E\{0\}) \cong \mathcal{M}(E)$ and $\pi_t(w) := w(t)$ denotes the projection on the endpoint of a path $w \in \mathcal{D}_E[0, t]$, then (a proof can be found in Section 2.2.3) the projection

$$(1.15) \qquad \mathcal{X}_t := \hat{\mathcal{X}}_t \circ \pi_t^{-1}, \qquad t \geq 0,$$

gives back the original $(G, \alpha, \beta)$-superprocess $\mathcal{X}$ started in $\mathcal{X}_0 = \mu$.

Likewise, if $X$ is a $(G, b, d)$-particle system in $E$ as defined in the last section, then the *historical binary branching particle system* $\hat{X}$ associated with $X$ is defined as the time-inhomogeneous binary branching particle system with time-dependent state space $\mathcal{N}(\mathcal{D}_E[0, t])$, with underlying motion $\hat{\xi}$, time-dependent branching rate $\hat{b}(t, w) := b(w(t))$ and time-dependent death rate $\hat{d}(t, w) := d(w(t))$. We call $\hat{X}$ the *historical $(G, b, d)$-particle system*. For a historical setting in the case of spatial Markov branching processes in discrete time, see, for instance, [13] or [18], Chapter 10. Viewed as a process in $\mathcal{N}(\mathcal{D}_E[0, \infty))$, $\hat{X}$ has cadlag sample paths. If $\hat{X}$ is started at time zero in $\hat{X}_0 = \nu \in \mathcal{N}(\mathcal{D}_E\{0\}) \cong \mathcal{N}(E)$, then the analogue of (1.15) gives back the (nonhistorical) $(G, b, d)$-particle system $X$ started in $X_0 = \nu$.

1.4. *Weighted superprocesses and compensated h-transforms*. We continue to assume that $\xi$ is a Feller process in a compact metrizable space $E$. Let $G$ be the generator of $\xi$, that is, $Gf := \lim_{t \to 0} t^{-1}(P_t f - f)$, where $P_t f(x) := E^x[f(\xi_t)]$ is the semigroup associated with $\xi$ and the domain $\mathcal{D}(G)$ of $G$ consists of all functions $f \in \mathcal{C}(E)$ for which the limit exists in the supremum norm. The following lemma, the proof of which can be found in Section 2.3.3, introduces *compensated h-transforms* of Feller processes.

LEMMA 3 (Compensated $h$-transform of a Feller process). *Let $G$ be the generator of a Feller process $\xi$ in a compact metrizable space $E$ and assume that $h \in \mathcal{D}(G)$ satisfies $h > 0$. Then the operator*

$$(1.16) \qquad G^h f := \frac{1}{h}(G(hf) - (Gh)f),$$

*with domain $\mathcal{D}(G^h) := \{f \in \mathcal{C}(E) : hf \in \mathcal{D}(G)\}$, is the generator of a Feller process $\xi^h$ on $E$. The laws of $\xi^h$ and $\xi$ are related by*

$$P^x[(\xi_s^h)_{s \in [0, t]} \in dw]$$

(1.17)
$$= \frac{h(w_t)}{h(x)} \exp\bigg(-\int_0^t \frac{Gh}{h}(w_s)\bigg) P^x[(\xi_s)_{s\in[0,t]} \in dw], \qquad t > 0, x \in E.$$

REMARK 4 ($h$-transforms). Doob's $h$-transform of a Feller process is the process with generator $\tilde{G}^h f := \frac{1}{h} G(hf)$ (see, e.g., [4], Section 2.VI.13, [24], formula (62.23) and [5], Section IX.4). Here $h$ is superharmonic (i.e., $Gh \leq 0$) and the $h$-transformed process has an additional local killing rate $Gh/h$. In our setup, it is natural to compensate for this killing by adding the term $-Gh/h$ in the definition of $G^h$. In this case, we can allow $h$ to be any positive function in the domain of $G$. A variant of the transformation in (1.16) appeared before in [15], Section 4. At least for diffusion processes, their transformation is equivalent to (1.16) if one chooses the logarithm of $h$ for their function $\xi$.

The following lemma, which was proved in a nonhistorical setting for superdiffusions in [7], describes the relation between weighted historical $(G, \alpha, \beta)$-superprocesses and compensated $h$-transforms.

LEMMA 5 (Weighted superprocess). *Let $\hat{\mathcal{X}}$ be the historical $(G, \alpha, \beta)$-superprocess and assume that $h \in \mathcal{D}(G)$, $h > 0$. Then the weighted process $\hat{\mathcal{X}}^h$, defined by*

(1.18) $$\hat{\mathcal{X}}_t^h(dw) := h(w_t)\hat{\mathcal{X}}_t(dw), \qquad t \geq 0,$$

*is the historical $(G^h, h\alpha, \beta + \frac{Gh}{h})$-superprocess.*

In particular, by formula (1.15), if $\mathcal{X}$ is the $(G, \alpha, \beta)$-superprocess, then $\mathcal{X}_t^h(dx) := h(x)\mathcal{X}_t(dx)$, $t \geq 0$, is the $(G^h, h\alpha, \beta + \frac{Gh}{h})$-superprocess. The proof of Lemma 5 is deferred to Section 2.3.4.

1.5. *Main results.* We are ready to state our first main result.

THEOREM 6 (Embedded particle system). *Let $E$ be a compact metrizable space, let $G$ be the generator of a Feller process in $E$ and $\alpha \in \mathcal{C}_+(E)$, $\beta \in \mathcal{C}(E)$. Assume that $h \in \mathcal{D}(G)$ satisfies $h > 0$ and, for some $\gamma \in \mathcal{C}_+(E)$,*

(1.19) $$Gh + \beta h - \alpha h^2 = -\gamma h.$$

*Then the historical $(G, \alpha, \beta)$-superprocess $\hat{\mathcal{X}}$ started in $\hat{\mathcal{X}}_0 = \mu \in \mathcal{M}(E)$ and the historical $(G^h, h\alpha, \gamma)$-particle system $\hat{X}$ started in $\hat{X}_0 = \text{Pois}(h\mu)$ can be coupled as processes such that*

(1.20) $$P[\hat{X}_t \in \cdot | (\hat{\mathcal{X}}_s)_{0 \leq s \leq t}] = P[\text{Pois}((h \circ \pi_t)\hat{\mathcal{X}}_t) \in \cdot | \hat{\mathcal{X}}_t] \qquad a.s. \ \forall t \geq 0.$$

It follows from (1.15) that the associated nonhistorical processes $\mathcal{X}$ and $X$ are related by (1.2). The phrase "coupled as processes" means that $(\hat{\mathcal{X}}_t)_{t\geq 0}$ and $(\hat{X}_t)_{t\geq 0}$ can be defined on the same probability space in such a way that (1.20) holds.

If $\hat{\mathcal{X}}$ and $\hat{X}$ are related by (1.20), then clearly the extinction of $\hat{\mathcal{X}}$ implies the extinction of $\hat{X}$ a.s. We now investigate when the converse conclusion can be drawn, that is, when $\hat{\mathcal{X}}$ and $\hat{X}$ can be coupled such that in addition to (1.20), eventually all mass of the superprocess $\mathcal{X}$ descends from particles in $X$. Set

$$(1.21) \qquad p(x) := -\log P^{\delta_x}[\mathcal{X}_t = 0 \ t\text{-eventually}], \qquad x \in E.$$

Here, $-\log 0 := \infty$ and we write $t$-*eventually* behind an event, depending on $t$, to denote the existence of a (random) time $\tau < \infty$ such that the event holds for all $t \geq \tau$. If no ambiguity is possible, we simply write eventually. It is not hard to check that $p$, defined by (1.21), satisfies (1.4). Therefore, we call $p$ the *infinitesimal survival probability* of $\mathcal{X}$. Note that

$$(1.22) \qquad P^{\delta_x}[\mathcal{X}_t = 0] = E^{\delta_x}[e^{-\langle \mathcal{X}_t, \infty \rangle}] = e^{-\mathcal{U}_t \infty(x)}, \qquad t \geq 0, x \in E.$$

The following proposition is proved in Section 3.1.3.

PROPOSITION 7 (Properties of the infinitesimal survival probability). *Consider $\mathcal{U} = \mathcal{U}(G, \alpha, \beta)$, where $G$, $\alpha$ and $\beta$ are as in Theorem 6, and let $p$ be given by (1.21). Assume that $\sup_{x \in E} \mathcal{U}_t \infty(x) < \infty$ for some $t > 0$. Then we have the following properties*:

(a) *Pointwise $\mathcal{U}_t \infty \downarrow p$ as $t \uparrow \infty$ and $\lim_{t \to \infty} \mathcal{U}_t f = p$ for all $f \in \mathcal{C}_+(E)$ with $f > 0$.*

(b) *For all $t \geq 0$, $\mathcal{U}_t p = p$.*

(c) *A function $f \in \mathcal{C}_+(E)$ satisfies $\mathcal{U}_t f = f$ for all $t \geq 0$ if and only if $f \in \mathcal{D}(G)$ and $f$ solves*

$$(1.23) \qquad\qquad\qquad Gf + \beta f - \alpha f^2 = 0.$$

(d) *If $\inf_{x \in E} p(x) > 0$, then $p$ is continuous and $p$ is the unique positive solution to* (1.23).

We now formulate our main theorem, which gives sufficient conditions for all mass of the superprocess $\mathcal{X}$ to descend eventually from particles in an embedded particle system $X$. We write $\pi_{[0,s]}$ to denote projection on $\mathcal{D}_E[0,s]$. By definition, the support $\mathrm{supp}(\mu)$ of a measure $\mu$ is the smallest closed set such that $\mu(\mathrm{supp}(\mu)^c) = 0$.

THEOREM 8 (Eventual descent from an embedded particle system). *Let $\hat{\mathcal{X}}$, $\hat{X}$ and $h$ be as in Theorem 6, and assume that $\mathcal{U} = \mathcal{U}(G, \alpha, \beta)$ satisfies* $\sup_{x \in E} \mathcal{U}_t \infty(x) < \infty$ for some $t > 0$. Then $p \leq h$. Moreover, $\hat{\mathcal{X}}$ and $\hat{X}$ may be coupled as processes such that (1.20) holds and such that, additionally,

$$(1.24) \quad \mathrm{supp}(\hat{X}_t) \supset \mathrm{supp}(\hat{\mathcal{X}}_r \circ \pi_{[0,t]}^{-1}) \qquad r\text{-eventually } \forall t \geq 0 \text{ a.s.}$$

If, moreover, $\inf_{x \in E} p(x) > 0$, then by Proposition 7 we may take $h = p$ in Theorem 6. In this case we have the following theorem:

THEOREM 9 (Trimmed tree of a superprocess). *Let $E$ be a compact metrizable space, let $G$ be the generator of a Feller process in $E$ and $\alpha \in \mathcal{C}_+(E)$, $\beta \in \mathcal{C}(E)$. Assume that $\mathcal{U} = \mathcal{U}(G, \alpha, \beta)$ satisfies $\sup_{x \in E} \mathcal{U}_t \infty(x) < \infty$ for some $t > 0$ and $\inf_{x \in E} p(x) > 0$. Then the historical $(G, \alpha, \beta)$-superprocess $\hat{\mathcal{X}}$ started in $\hat{\mathcal{X}}_0 = \mu \in \mathcal{M}(E)$ and the historical $(G^p, p\alpha, 0)$-particle system $\hat{X}$ started in $\hat{X}_0 = \mathrm{Pois}(p\mu)$ can be coupled as processes such that*

$$(1.25) \quad P[\hat{X}_t \in \cdot | (\hat{\mathcal{X}}_s)_{0 \leq s \leq t}] = P[\mathrm{Pois}((p \circ \pi_t)\hat{\mathcal{X}}_t) \in \cdot | \hat{\mathcal{X}}_t] \qquad \text{a.s. } \forall t \geq 0$$

*and*

$$(1.26) \quad \mathrm{supp}(\hat{X}_t) = \mathrm{supp}(\hat{\mathcal{X}}_r \circ \pi_{[0,t]}^{-1}), \qquad r\text{-eventually } \forall t \geq 0 \text{ a.s.}$$

If $\hat{\mathcal{X}}$ and $\hat{X}$ are coupled as in Theorem 9, then we say that $\hat{X}$ is the *trimmed tree* of $\hat{\mathcal{X}}$. If $\mathcal{X}_t = \hat{\mathcal{X}}_t \circ \pi_t^{-1}$ and $X_t = \hat{X}_t \circ \pi_t^{-1}$ are the associated nonhistorical processes, then we also call $X$ the trimmed tree of $\mathcal{X}$. Note that the death rate of $X$ is zero, that is, $X$ is a binary splitting particle system.

REMARK 10 (Checking the assumptions on $\mathcal{U}_t \infty$ and $p$). Upper bounds on $\mathcal{U}_t \infty$ and lower bounds on $p$ can be found, in practical situations, by finding solutions to an appropriate differential inequality; see Lemmas 12 and 25.

1.6. *Finite ancestry.* In this section, we investigate the assumption in Theorems 8 and 9 that $\sup_{x \in E} \mathcal{U}_t \infty(x) < \infty$ for some $t > 0$. In particular, we show that this assumption is equivalent to the statement that all mass of the superprocess $\mathcal{X}$ descends eventually from finitely many ancestors, in some sense.

To do this, we need to equip the historical $(G, \alpha, \beta)$-superprocess $\hat{\mathcal{X}}$ with some additional structure that makes it possible to distinguish different ancestors. To this aim, set $E' := E \times [0, 1]$. Define a Feller process $\xi' = (\xi, \eta)$ on $E'$, where for given initial conditions $(x, y) \in E \times [0, 1]$, $\xi$ is the Feller process with generator $G$ started in $x$, and $\eta_t := y$, $t \geq 0$. Put $\alpha'(x, y) := \alpha(x)$ and $\beta'(x, y) := \beta(x)$. Let $\hat{\mathcal{X}}'$ denote the historical $(G', \alpha', \beta')$-superprocess. Then the formula

$$(1.27) \qquad \hat{\mathcal{X}}_t := \hat{\mathcal{X}}'_t \circ \psi_t^{-1}, \qquad t \geq 0,$$

gives back the original historical $(G,\alpha,\beta)$-superprocess $\hat{\mathcal{X}}$, where $\psi_t$ denotes the projection from $\mathcal{D}_{E\times[0,1]}[0,t]$ to $\mathcal{D}_E[0,t]$. The following lemma is proved in Section 3.2.3. Here $\pi_0(w) := w(0)$ denotes the projection on the starting point of a path $w$ in $\mathcal{D}_E[0,t]$ or $\mathcal{D}_{E'}[0,t]$.

LEMMA 11 (Finite ancestry). *Let $\hat{\mathcal{X}}$ be the historical $(G,\alpha,\beta)$-superprocess, let $\hat{\mathcal{X}}'$ be the extended historical $(G',\alpha',\beta')$-superprocess just defined and $\mathcal{U} = \mathcal{U}(G,\alpha,\beta)$. Let $\ell$ denote Lebesgue measure on $[0,1]$. Then we have the relations* (i) $\Leftrightarrow$ (ii) $\Rightarrow$ (iii), *where*

$$(1.28) \quad \begin{aligned} &\text{(i)} \sup_{x\in E} \mathcal{U}_t\infty(x) < \infty \qquad \text{for some } t > 0, \\ &\text{(ii)} \ P^{0,\mu\otimes\ell}[\mathrm{supp}(\hat{\mathcal{X}}'_t \circ \pi_0^{-1}) \text{ is finite eventually}] = 1 \qquad \forall \mu \in \mathcal{M}(E), \\ &\text{(iii)} \quad P^{0,\mu}[\mathrm{supp}(\hat{\mathcal{X}}_t \circ \pi_0^{-1}) \text{ is finite eventually}] = 1 \qquad \forall \mu \in \mathcal{M}(E). \end{aligned}$$

We interpret $\mathrm{supp}(\hat{\mathcal{X}}'_t \circ \pi_0^{-1})$ as the ancestors at time 0 of the population of $\mathcal{X}$ at time $t$. We have extended the underlying space $E$ to make sure that different ancestors live a.s. on different positions. Note that if $E$ is finite, then (iii) is always trivially fulfilled even when (i) fails.

For many superprocesses, it is actually the case that

$$(1.29) \qquad \sup_{x\in E} \mathcal{U}_t\infty(x) < \infty \qquad \forall t > 0.$$

A sufficient, but not necessary condition for (1.29) is that $\alpha$ is bounded away from zero. The sufficiency follows from the following bound (see, e.g., [14], Lemma 11).

LEMMA 12 (Extinction estimate). *Set $\underline{\alpha} := \inf_{x\in E}\alpha(x)$ and $\overline{\beta} := \sup_{x\in E}\beta(x)$. If $\underline{\alpha} > 0$, then*

$$(1.30) \quad \mathcal{U}_t\infty \leq \frac{\overline{\beta}}{\underline{\alpha}(1-e^{-\overline{\beta}t})}, \qquad \overline{\beta} \neq 0 \qquad \text{and} \qquad \mathcal{U}_t\infty \leq \frac{1}{\underline{\alpha}t}, \qquad \overline{\beta} = 0.$$

On the other hand, it is possible for a $(G,\alpha,\beta)$-superprocess to satisfy (1.29) while $\underline{\alpha} = 0$ (see [14], Lemmas 5 and 6).

The following consequence of (1.29) is proved in Section 3.2.3.

LEMMA 13 (Finite ancestry and preserved past property). *If $\mathcal{X}$ satisfies* (1.29), *then*

$$(1.31) \quad \begin{aligned} &\text{(i) } \mathrm{supp}\,(\hat{\mathcal{X}}_r \circ \pi_{[0,t]}^{-1}) \text{ is finite} \qquad \forall 0 \leq t < r \text{ a.s.}, \\ &\text{(ii) } \mathrm{supp}\,(\hat{\mathcal{X}}_r \circ \pi_{[0,t]}^{-1}) \supset \mathrm{supp}\,(\hat{\mathcal{X}}_{r'} \circ \pi_{[0,t]}^{-1}) \qquad \forall 0 \leq t < r \leq r' \text{ a.s.} \end{aligned}$$

In view of Lemma 11 and (1.31)(i) we say that a superprocess $\mathcal{X}$ has the *finite ancestry property* if $\mathcal{X}$ satisfies (1.29). Note that (1.31)(ii) says that lines of descent (up to a given time $s$) can become extinct, but no new ones are created. This statement may seem obvious, but some care is needed regarding the order of the $\forall$ and the a.s. in the statements. In (1.31)(ii), we claim that the same zero set works for all times $t, r, r'$ such that $0 \leq t < r \leq r'$. One cannot simply argue by continuity here, because the support of a measure $\mu$ is not a continuous function of $\mu$. Note that if the superprocess $\mathcal{X}$ in Theorem 9 has the finite ancestry property, then a.s. the sets $\mathrm{supp}(\hat{\mathcal{X}}_r \circ \pi_{[0,t]}^{-1})$ in (1.26) are finite for all $r > t$ and decrease to $\mathrm{supp}(\hat{X}_t)$ as $r \uparrow \infty$.

1.7. *Methods, discussion and outline of the proofs.* Our results have obvious applications in the study of (local) extinction and survival of superprocesses. For superdiffusions, extinction properties were studied by Engländer and Pinsky [7]. Parallel to the present paper, Engländer and Kyprianou [6] investigated local survival and local exponential growth of superdiffusions. The first paper uses more analytic tools, while the second is more probabilistic in nature.

While our methods are more probabilistic, some of our results are close in spirit to the work in [7]. As we already mentioned, the weighted superprocess $(p\mathcal{X}_t)_{t \geq 0}$ with $p$ as in (1.4) plays an important role in [7]. Also, their Theorem 4.4(a) is not surprising in view of our Theorem 9, although their setup and ours do not completely overlap. Their Theorem 3.1 describes properties of the function $p$ similar (but not identical) to our Proposition 7. Since our underlying motion is a general Feller process which does not have the good smoothing properties of uniformly elliptic diffusions, we have to be more careful about the sense in which $p$ solves equation (1.23).

The main tool in [6] is an expression [their Theorem 5(i)] that says (in the language of log-Laplace functionals) that a certain change of measure of a superdiffusion yields back the same superdiffusion with an additional immigration term coming from a single particle. In their introduction, Engländer and Kyprianou discussed the possible use of Poissonization relations for their aims, but rejected them on the ground that relation between the laws of $\mathcal{X}_t$ and $X_t$ at fixed times $t$ are not enough to relate the long-time behavior of $\mathcal{X}$ and $X$. A central aim of our work is to overcome such shortcomings of the usual Poissonization formulas. Another aim, of course, is to allow more general underlying motions than diffusions.

The main ideas behind our proofs of Theorems 6, 8 and 9 are the simple observations about Poissonization and weighting of superprocesses in Lemmas 1 and 5, respectively. Our strategy is to construct a version of the superprocess with so much additional structure that we can distinguish all ancestors of the population alive at a given time. For such a sufficiently enriched process, we then explicitly identify the trimmed tree and check that it is a binary splitting particle system. This is done in Proposition 39 and Lemma 40 in Sections 3.3.1 and 3.3.2, respectively. The essential step, where a coupling of $\mathcal{X}_t$ and $X_t$ for fixed $t$ is improved to a coupling of $\mathcal{X}$ and $X$ as processes, occurs in the proof of Lemma 40. Forgetting step by step some of the added structure, we then arrive at Theorems 6, 8 and 9.

Interesting side results of this approach are a number of lemmas about the lines of descent of a superprocess, notably Lemma 13, which may seem intuitively obvious, but to our knowledge has not been proved before. On the other hand, our approach does not make any statements about the transition probabilities of the joint process $(\mathcal{X}_t, X_t)_{t \geq 0}$, when $\mathcal{X}$ and $X$ (and their historical counterparts) are coupled as in Theorem 6. Another possible approach to our Theorem 6 (not followed in this paper) would be to specify a joint Markov evolution for $(\mathcal{X}, X)$ and then show that if the process is started in a state such that $X_0 = \text{Pois}(h\mathcal{X}_0)$, then $X_t = \text{Pois}(h\mathcal{X}_t)$ for all $t \geq 0$. Here, $X$ would be an autonomous binary branching particle system, while $\mathcal{X}$ would be a superprocess with an additional mass creation on the positions of the particles in $X$.

Our results can be generalized in several directions. If the space $E$ is not compact but *locally compact,* then generalizations of our results can be derived using the compactification technique from Remark 2. This requires, however, that the functions $h$ in Theorems 6 and 8 or $p$ in Theorem 9 are uniformly bounded away from zero, and hence can be extended to positive continuous functions on some compactification of $E$. Truly local versions of our results, where $h$ and $p$ are only required to be locally bounded away from zero, are somewhat more subtle. We hope to handle these in a forthcoming paper.

A lot of our proofs work for superprocesses whose underlying motion is a general Hunt process on a Polish space, and whose activity and growth parameter are bounded and measurable, but we do not know how to treat compensated $h$-transforms and weighted superprocesses (Lemmas 3 and 5) in this context.

The proofs are organized as follows. After settling some notational and topological issues in Section 2.1, we introduce formally historical processes in Section 2.2 and collect some of their elementary properties. Section 2.3 treats compensated $h$-transforms and weighted superprocesses. Section 3.1 is devoted to the infinitesimal survival probability $p$. Section 3.2 collects some basic facts about surviving lines of descent. In Section 3.3, finally, we prove our main results.

**2. Prerequisites on superprocesses.**

2.1. *Topological preliminaries.* Let $E$ be a Polish space (i.e., $E$ is a separable topological space and there exists a complete metric generating the topology). We always equip $E$ with the Borel $\sigma$-field $\mathcal{B}(E)$. We let $B(E)$, $B_+(E)$ and $B_{[0,1]}(E)$ denote the spaces of bounded, bounded nonnegative and $[0,1]$-valued, real measurable functions on $E$, respectively. If a countable collection of functions $\{f_i : i \geq 1\} \subset B(E)$ separates points, then $\mathcal{B}(E) = \sigma(f_i : i \geq 1)$ (see [23], Lemma II.18). We remind the reader of the fact that a subspace $F$ of a Polish space $E$ is itself Polish in the induced topology if and only if $F$ is a $G_\delta$-subset of $E$, that is, a countable intersection of open sets ([1], Section 6, Theorem 1).

Let $\mathcal{C}_b(E)$ denote the space of bounded continuous real functions on $E$. We write $\mathcal{M}(E)$ for the space of finite measures on $E$, equipped with the topology of weak convergence (with weak convergence denoted by $\Rightarrow$), under which $\mathcal{M}(E)$ is a Polish space ([9], Theorem 3.1.7). Recall that by definition $\mu_n \Rightarrow \mu$ if and only if $\langle \mu_n, f \rangle \to \langle \mu, f \rangle$ for all $f \in \mathcal{C}_b(E)$. Note that the topology on $\mathcal{M}(E)$ does not depend on the choice of the metric on $E$. The Borel $\sigma$-field on $\mathcal{M}(E)$ is generated by the mappings $\mu \mapsto \mu(A)$, $A \in \mathcal{B}(E)$ (cf. [21], Lemma 3.2.3). If $F \subset E$ is measurable, we identify $\mathcal{M}(F)$ with the space $\{\mu \in \mathcal{M}(E) : \mu(E \setminus F) = 0\}$. In particular, when $F$ is a $G_\delta$-subset of $E$ (and therefore Polish in the induced topology), then the topology of weak convergence on $\mathcal{M}(F)$ coincides with the induced topology from its embedding in $\mathcal{M}(E)$. By $\mathcal{M}_1(E) \subset \mathcal{M}(E)$ we denote the space of probability measures; $\mathcal{N}(E) \subset \mathcal{M}(E)$ denotes the space of finite point measures on $E$.

We denote by $\mathcal{D}_E[0, \infty)$ the space of cadlag (i.e., right-continuous with existing left limits) functions $w : [0, \infty) \to E$, equipped with the Skorohod topology. This is the $J_1$ topology defined in [25]. The space $\mathcal{D}_E[0, \infty)$ is Polish ([9], Theorem 3.5.6). One has $w_n \to w$ in $\mathcal{D}_E[0, \infty)$ if and only if for each $T > 0$ there exists a sequence of strictly increasing, continuous $\lambda_n : [0, T] \to [0, \infty)$ with $\lambda_n(0) = 0$, such that

$$\lim_{n \to \infty} \sup_{t \in [0,T]} |\lambda_n(t) - t| = 0 \tag{2.1}$$

and such that (cf. [9], Proposition 3.5.3)

$$w_n(\lambda_n(t_n)) \to \begin{cases} w(t), & \text{whenever } t_n \downarrow t, \\ w(t-), & \text{whenever } t_n \uparrow t, \end{cases} \quad t_n, t \in [0, T]. \tag{2.2}$$

Note that the topology on $\mathcal{D}_E[0, \infty)$ does not depend on the choice of the metric on $E$.

2.2. *Historical processes.*

2.2.1. *Hunt processes.* Let $E$ be a Polish space and let $(P_t)_{t\geq 0}$ be a measurable transition probability on $E$. That is, $(t,x) \mapsto P_t(x,\cdot)$ is a (Borel) measurable map from $[0,\infty) \times E$ into $\mathcal{M}_1(E)$, $P_0(x,\cdot) = \delta_0$ for all $x \in E$ and the operators

$$(2.3) \qquad P_t f(x) := \int_E P_t(x,dy) f(y), \qquad t \geq 0, x \in E, f \in B(E),$$

form a semigroup: $P_t P_s f = P_{t+s} f$ for all $s,t \geq 0$, $f \in B(E)$.

Assume that $(P_t)_{t\geq 0}$ is the transition probability (equivalently the semigroup) of a Markov process with cadlag sample paths in $E$, that is, for every $x \in E$ there exists a $\mathcal{D}_E[0,\infty)$-valued random variable $\xi^x$, unique in distribution, such that $\xi_0^x = x$ and

$$(2.4) \qquad E[f(\xi_t^x)|\mathcal{F}_s] = (P_{t-s}f)(\xi_s^x) \qquad \text{a.s., } 0 \leq s \leq t, f \in B(E),$$

where $(\mathcal{F}_t)_{t\geq 0}$ denotes the filtration generated by $\xi^x$. By definition, the Markov process with transition probability $(P_t)_{t\geq 0}$ is a *Hunt process* if, for every $x \in E$, the following statements hold (see [24], Theorem I.7.4 and Definition V.47.3):

(2.5)
(i) *Right property.* For every $t > 0$ and $f \in B(E)$,
  the map $[0,t) \ni s \mapsto P_{t-s}f(\xi_s^x)$ is a.s. right-continuous.
(ii) *Quasi left-continuity.* For every increasing sequence of
  $\mathcal{F}_{\cdot+}$ stopping times $\tau_n \uparrow \tau$, we have $\xi_{\tau_n}^x \to \xi_\tau^x$ a.s. on $\{\tau < \infty\}$.

Here $\mathcal{F}_{\cdot+} = (\mathcal{F}_{t+})_{t\geq 0}$ denotes the right-continuous modification of $(\mathcal{F}_t)_{t\geq 0}$. The right property implies the strong Markov property ([24], Theorem I.7.4). Conditions (2.5)(i) and (2.5)(ii) are properties of the law $P^x := \mathcal{L}(\xi^x)$ of $\xi^x$ only and, therefore, being a Hunt process is a property of the transition probability. It suffices to check (2.5)(i) for all $f \in \mathcal{C}_b(E)$ ([24], Theorem I.7.4). We identify a Hunt process with the collection of probability measures $(P^x)^{x \in E}$. If $(\xi^x)^{x \in E}$ is a collection of $\mathcal{D}_E[0,\infty)$-valued random variables with laws $\mathcal{L}(\xi^x) = P^x$, $x \in E$, then with a slight abuse of terminology we say that $(\xi^x)^{x \in E}$ is a Hunt process (regardless of a possible dependence structure between the $\xi^x$).

We also need time-inhomogeneous Hunt processes with a time-dependent state space $E_t$. We assume that the $E_t$ are (or can be identified with) subsets of some Polish space $E$ and that the set $\dot E := \{(t,x) \in [0,\infty) \times E : x \in E_t\}$ is a $G_\delta$-subset of $[0,\infty) \times E$ (and therefore Polish in the induced topology). Let $W_{[s,\infty)} := \{w \in \mathcal{D}_E[s,\infty) : w_t \in E_t\ \forall t \geq s\}$ denote the space of all possible paths the process can follow after time $s$. Generalizing our previous definition, we say that a collection of random variables $(\xi^{s,x})^{(s,x)\in\dot E}$, where $\xi^{s,x}$ takes values in $W_{[s,\infty)}$, is a *time-inhomogeneous Hunt process* if the collection of random variables $(\dot\xi^{(s,x)})^{(s,x)\in\dot E}$ defined by

$$(2.6) \qquad \dot\xi_t^{(s,x)} := (s+t, \xi_{s+t}^{s,x}), \qquad (s,x) \in \dot E, t \geq 0,$$

is a (time-homogeneous) Hunt process in $\dot E$. If $(\xi^{s,x})^{(s,x)\in \dot E}$ is a time-inhomogeneous Hunt process, then we write $P_{s,t}(x,\cdot):=P[\xi_t^{s,x}\in\cdot]$ and we let $P_{s,t}\colon B(E_t)\to B(E_s)$ denote the operator

$$(2.7) \qquad P_{s,t}f(x):=\int_{E_t} P_{s,t}(x,dy)f(y), \qquad x\in E_s, f\in B(E_t).$$

By a slight abuse of terminology, we call $(P_{s,t})_{t\geq s\geq 0}$ the (time-inhomogeneous) semigroup associated with $(\xi^{s,x})^{(s,x)\in \dot E}$. (Such time-inhomogeneous semigroups are sometimes called transition functions.)

2.2.2. *Superprocesses with Hunt underlying motion.* Let $\xi$ be a (time-homogeneous) Hunt process in a Polish space $E$ with semigroup $(P_t)_{t\geq 0}$ and assume that $\alpha\in B_+(E)$, $\beta\in B(E)$. Then, for every $f\in B_+(E)$, there exists a unique $\mathcal{B}([0,\infty)\times E)$-measurable nonnegative function $u$ which is bounded on $[0,T]\times E$ for all $T>0$, solving the Cauchy integral equation

$$(2.8) \qquad u_t = P_t f + \int_0^t P_{t-s}(\beta u_s - \alpha u_s^2)\, ds, \qquad t\geq 0$$

([10], Proposition 2.3). Moreover, it was shown ([10], Corollary 3.6) that there exists a unique (in law) Hunt process $(\mathcal{X}^\mu)^{\mu\in\mathcal{M}(E)}$, with continuous sample paths, such that

$$(2.9) \qquad E^\mu[e^{-\langle \mathcal{X}_t, f\rangle}] = e^{-\langle \mu, \mathcal{U}_t f\rangle}, \qquad t\geq 0, \mu\in\mathcal{M}(E), f\in B_+(E),$$

where $\mathcal{U}_t f := u_t$, $t\geq 0$, and $u$ solves (2.8). We call $\mathcal{X}$ the superprocess with underlying motion $\xi$, activity $\alpha$ and growth parameter $\beta$, or, for short, the $(\xi,\alpha,\beta)$-*superprocess*, and we call $\mathcal{U}=\mathcal{U}(\xi,\alpha,\beta)$ its *log-Laplace semigroup*. By monotone convergence, $\mathcal{U}_t f$ can be defined unambiguously such that (2.9) holds for any measurable $f\colon E\to[0,\infty]$ ([14], Lemma 9).

We list some elementary properties of $(\xi,\alpha,\beta)$-superprocesses that we need later. The following lemma is an easy consequence of (2.9).

LEMMA 14 (Branching property). *Let $\mu_1,\mu_2\in\mathcal{M}(E)$, and let $\mathcal{X}^{\mu_1}$ and $\mathcal{X}^{\mu_2}$ be independent copies of the $(\xi,\alpha,\beta)$-superprocess started in $\mu_1$ and $\mu_2$, respectively. Then*

$$(2.10) \qquad \mathcal{X}_t^{\mu_1+\mu_2} := \mathcal{X}_t^{\mu_1} + \mathcal{X}_t^{\mu_2}, \qquad t\geq 0,$$

*is the $(\xi,\alpha,\beta)$ superprocess started in $\mu_1+\mu_2$.*

The following lemma was proved in [10], Proposition 2.7.

LEMMA 15 (Moment formulas). *For every $f \in B(E)$, there exists a unique $\mathcal{B}([0, \infty) \times E)$-measurable function $v$ which is bounded on $[0, T] \times E$ for all $T > 0$, such that*

$$(2.11) \qquad v_t = P_t f + \int_0^t P_{t-s}(\beta v_s) \, ds, \qquad t \geq 0.$$

*The formula $\mathcal{V}_t f := v_t$ defines a (linear) semigroup $(\mathcal{V}_t)_{t \geq 0}$ on $B(E)$. We have*

$$(2.12) \quad \mathcal{V}_t f(x) = E^x \bigg[ f(\xi_t) \exp \bigg( \int_0^t \beta(\xi_s) \, ds \bigg) \bigg], \qquad t \geq 0, x \in E, f \in B(E).$$

*Moreover, for all $t \geq 0, f, g \in B(E)$,*

$$(2.13) \quad \begin{aligned} &\text{(i) } E^\mu[\langle \mathcal{X}_t, f \rangle] = \langle \mu, \mathcal{V}_t f \rangle, \\ &\text{(ii) } \mathrm{Cov}^\mu(\langle \mathcal{X}_t, f \rangle, \langle \mathcal{X}_t, g \rangle) = 2 \int_0^t ds \langle \mu, \mathcal{V}_s(\alpha(\mathcal{V}_{t-s} f)(\mathcal{V}_{t-s} g)) \rangle. \end{aligned}$$

The following lemma is an easy consequence of Lemma 15 and the fact that $0 \leq \mathcal{V}_t f \leq e^{\|\beta\| t} \|P_t f\|$ for all $f \in B_+(E)$ (where $\|\cdot\|$ denotes the supremum norm).

LEMMA 16 (Absolute continuity of moment measures). *Let $\mu$ be a probability measure on $E$ and $m \geq 0$. Then, for $t \geq 0$,*

$$(2.14) \quad \begin{aligned} &\text{(i) } E^{m\mu}[\mathcal{X}_t] \ll P^\mu[\xi_t \in \cdot], \\ &\text{(ii) } E^{m\mu}[\mathcal{X}_t \otimes \mathcal{X}_t] \ll P^\mu[\xi_t \in \cdot] \otimes P^\mu[\xi_t \in \cdot] + Q_t^\mu, \end{aligned}$$

*where $Q_t^\mu$ is the measure on $E \times E$ defined as*

$$(2.15) \qquad Q_t^\mu := \int_0^t ds \int_E P^\mu[\xi_s \in dx](P^x[\xi_{t-s} \in \cdot] \otimes P^x[\xi_{t-s} \in \cdot]).$$

A measure $\mu \in \mathcal{M}(E)$ is atomless [i.e., $\mu(\{x\}) = 0$ for all $x \in E$] if and only if

$$(2.16) \qquad \mu \otimes \mu(\{(x_1, x_2) \in E \times E : x_1 = x_2\}) = 0.$$

The following lemma follows from formulas (2.14)(ii) and (2.16).

LEMMA 17 (Atomless superprocess). *Assume that $P^x[\xi_t \in \cdot]$ is atomless for every $t > 0$ and $x \in E$. Then $\mathcal{X}_t$ is atomless a.s. for every $t > 0$ and initial state $\mu \in \mathcal{M}(E)$.*

Our next lemma is the following:

LEMMA 18 (Image property). *Let $E, F$ be Polish spaces, let $\psi : E \to F$ be continuous and let $\xi = (\xi^x)^{x \in E}$ and $\eta = (\eta^y)^{y \in F}$ be Hunt processes in $E$ and $F$, respectively, satisfying*

(2.17) $$\psi(\xi_t^x) = \eta_t^{\psi(x)}, \qquad x \in E, t \geq 0.$$

*Assume that $\alpha_F \in B_+(F)$ and $\beta_F \in B(F)$, and let $\alpha_E \in B_+(E)$ and $\beta_E \in B(E)$ be given by*

(2.18) $$\alpha_E := \alpha_F \circ \psi \quad \text{and} \quad \beta_E := \beta_F \circ \psi.$$

*Let $\mathcal{X}$ be the $(\xi, \alpha_E, \beta_E)$-superprocess with initial state $\mu \in \mathcal{M}(E)$. Then*

(2.19) $$\mathcal{Y}_t := \mathcal{X}_t \circ \psi^{-1}, \qquad t \geq 0,$$

*is the $(\eta, \alpha_F, \beta_F)$-superprocess with initial state $\mu \circ \psi^{-1}$.*

PROOF. Let $P^E$ and $P^F$ denote the semigroups associated with the processes $\xi$ and $\eta$, respectively. Formula (2.17) implies that $P_t^E(f \circ \psi) = (P_t^F f) \circ \psi$ for all $f \in B(F)$. Using this fact and (2.18), it is not hard to show that also $\mathcal{U}_t^E(f \circ \psi) = (\mathcal{U}_t^F f) \circ \psi$ for all $f \in B_+(F)$, where $\mathcal{U}^E = \mathcal{U}(\xi, \alpha_E, \beta_E)$ and $\mathcal{U}^F = \mathcal{U}(\eta, \alpha_F, \beta_F)$ are the log-Laplace semigroups of $\mathcal{X}$ and $\mathcal{Y}$, respectively. Let $(\mathcal{F}_t)_{t \geq 0}$ be the filtration generated by $\mathcal{X}$. Then, for all $0 \leq s \leq t$,

(2.20)
$$\begin{aligned}
E[\exp(-\langle \mathcal{X}_t \circ \psi^{-1}, f \rangle) | \mathcal{F}_s] \\
= E[\exp(-\langle \mathcal{X}_t, f \circ \psi \rangle) | \mathcal{F}_s] = \exp\big(-\langle \mathcal{X}_s, U_{t-s}^E(f \circ \psi) \rangle\big) \\
= \exp\big(-\langle \mathcal{X}_s, (U_{t-s}^F f) \circ \psi \rangle\big) \\
= \exp(-\langle \mathcal{X}_s \circ \psi^{-1}, U_{t-s}^F f \rangle), \qquad f \in B_+(F).
\end{aligned}$$

This shows that $(\mathcal{X}_t \circ \psi^{-1})_{t \geq 0}$ is a Markov process and that its transition probabilities coincide with those of the $(\eta, \alpha_F, \beta_F)$-superprocess. Since $\psi$ is continuous, $\mathcal{X}_t \circ \psi^{-1}$ has continuous sample paths. □

The following simple observation will be useful later.

LEMMA 19 (Preserved sets). *Let $\mathcal{X}$ be the $(\xi, \alpha, \beta)$-superprocess.*

(a) *If $F \subset E$ is measurable and $P^x[\xi_t \in F] = 1 \, \forall t \geq 0 \, (x \in F)$, then*

(2.21) $$P^\mu[\mathcal{X}_t \in \mathcal{M}(F)] = 1 \qquad \forall t \geq 0, \mu \in \mathcal{M}(F).$$

(b) *If $F \subset E$ is a $G_\delta$-set and $P^x[\xi_t \in F \, \forall t \geq 0, \xi_{t-} \in F \, \forall t > 0] = 1$, $x \in F$, then*

(2.22) $$P^\mu[\mathcal{X}_t \in \mathcal{M}(F) \, \forall t \geq 0] = 1, \qquad \mu \in \mathcal{M}(F).$$

PROOF. Statement (a) follows from (2.14)(i), while (b) follows by applying Lemma 18 to the inclusion map $F \subset E$, where we use that the restriction of $\xi$ to $F$ is again a Hunt process. The assumption that $F$ is a $G_\delta$-set guarantees that $F$ is a Polish space and that the event $\{\mathcal{X}_t \in \mathcal{M}(F) \ \forall t \geq 0\}$ is Borel measurable. $\square$

We conclude this section by constructing superprocesses with time-inhomogeneous underlying motion. Let $\xi = (\xi^{s,x})^{(s,x) \in \dot{E}}$ be a time-inhomogeneous Hunt process as defined at the end of the last section, and assume that $\dot{\alpha} \in B_+(\dot{E})$ and $\dot{\beta} \in B(\dot{E})$. Let $\dot{\xi}$ be the time-homogeneous Hunt process in (2.6) and let $\dot{\mathcal{X}}$ denote the $(\dot{\xi}, \dot{\alpha}, \dot{\beta})$ superprocess. Using Lemma 16 we see that $\dot{\mathcal{X}}_t^{\delta_s \otimes \mu}$ is concentrated on $\{s+t\} \times E_{s+t}$ a.s. $\forall t \geq 0$. Since $\dot{\mathcal{X}}^{\delta_s \otimes \mu}$ has continuous sample paths and since $\{\delta_t \otimes \mu : t \geq 0, \mu \in \mathcal{M}(E_t)\} \subset \mathcal{M}(\dot{E})$ is closed, there exists a process $\mathcal{X}^{s,\mu}$ with continuous sample paths in $\mathcal{M}(E)$ such that $\mathcal{X}_{s+t}^{s,\mu} \in \mathcal{M}(E_{s+t})$ for all $t \geq 0$ and

$$\dot{\mathcal{X}}_t^{\delta_s \otimes \mu} = \delta_{s+t} \otimes \mathcal{X}_{s+t}^{s,\mu}. \tag{2.23}$$

Set $\dot{\mathcal{M}} := \{(t, \mu) \in [0, \infty) \times \mathcal{M}(E) : \mu \in \mathcal{M}(E_t)\}$. It is not hard to check that $\mathcal{X} = (\mathcal{X}^{s,\mu})^{(s,\mu) \in \dot{\mathcal{M}}}$ is a time-inhomogeneous Hunt process with continuous sample paths, and

$$E^{s,\mu}[e^{-\langle \mathcal{X}_t, f \rangle}] = e^{-\langle \mu, \mathcal{U}_{s,t} f \rangle}, \qquad t \geq s \geq 0, \mu \in \mathcal{M}(E_s), f \in B_+(E_t), \tag{2.24}$$

where $(\mathcal{U}_{s,t} f)_{s \in [0,t]} =: u \in B_+(\{(s,x) \in [0,t] \times E : x \in E_s\})$ solves the equation

$$u_s = P_{s,t} f + \int_s^t P_{s,r}(\beta_r u_r - \alpha_r u_r^2) \, dr, \qquad s \in [0, t]. \tag{2.25}$$

Here $\alpha_t(x) := \dot{\alpha}(t, x)$, $\beta_t(x) := \dot{\beta}(t, x)$ $((t, x) \in \dot{E})$ and $(P_{s,t})_{t \geq s \geq 0}$ is the (time-inhomogeneous) semigroup associated with $\xi$. We call $\mathcal{X}$ the (time-inhomogeneous) $(\xi, \alpha_t, \beta_t)$-superprocess and call $(\mathcal{U}_{s,t})_{t \geq s \geq 0}$ the (time-inhomogeneous) *log-Laplace semigroup* associated with $\mathcal{X}$.

2.2.3. *Historical superprocesses.* Let $\xi = (\xi^x)^{x \in E}$ be a Hunt process in a Polish space $E$ and let $\hat{\xi} = (\hat{\xi}^{s,w})^{s \geq 0, w \in \mathcal{D}_E[0,s]}$ be the associated path process, defined as in (1.14). Identify, as usual, $\mathcal{D}_E[0, s]$ with the subspace of $\mathcal{D}_E[0, \infty)$ consisting of paths stopped at time $s$ and define $\tilde{E} \subset [0, \infty) \times \mathcal{D}_E[0, \infty)$ by

$$\tilde{E} := \{(s, w) : s \geq 0, w \in \mathcal{D}_E[0, s]\}. \tag{2.26}$$

Then $(\hat{\xi}^{s,w})^{(s,w) \in \tilde{E}}$ is a time-inhomogeneous Hunt process (see [3], Proposition 2.1.2). If $\mathcal{X}$ is a $(\xi, \alpha, \beta)$-superprocess, then by definition the *historical $(\xi, \alpha, \beta)$-superprocess* $\hat{\mathcal{X}}$ is the (time-inhomogeneous) $(\hat{\xi}, \hat{\alpha}_t, \hat{\beta}_t)$-superprocess, where $\hat{\alpha}_t(w) := \alpha(w(t))$ and $\hat{\beta}_t(w) := \beta(w(t)), (t, w) \in \tilde{E}$. We are now in a situation where we can prove some of the elementary properties of historical superprocesses mentioned in Section 1.

PROOF OF (1.15). If $\hat{\xi}$ is the path process associated with a Hunt process $\xi$, started at time $s \geq 0$ in $w \in \mathcal{D}_E[0,s]$, then $\xi_t := \pi_{s+t}(\hat{\xi}_{s+t})$, $t \geq 0$, gives back the original Hunt process $\xi$ started in $\pi_s(\hat{\xi}_s)$. Moreover, the map $(t,w) \mapsto w(t)$ from $\tilde{E}$ into $E$ is continuous. (Note that this is true even though the map $w \mapsto w(t)$ from $\mathcal{D}_E[0,\infty)$ into $E$ is in general discontinuous.) Therefore, Lemma 18 (the image property of superprocesses) shows that if $(\hat{\mathcal{X}}_t)_{t \geq s}$ is the historical $(\xi, \alpha, \beta)$-superprocess started at time $s \geq 0$ in $\hat{\mu} \in \mathcal{D}_E[0,s]$, then

$$\mathcal{X}_t := \hat{\mathcal{X}}_{s+t} \circ \pi_{s+t}^{-1}, \qquad t \geq 0, \tag{2.27}$$

is the (nonhistorical) $(\xi, \alpha, \beta)$-superprocess started in $\hat{\mu} \circ \pi_s^{-1}$. □

One of the driving ideas behind the development of historical superprocesses has been the desire to have a means to distinguish those parts of the population that descend from different ancestors. However, all that a path in $\mathcal{D}_E[0,t]$ tells us is where in space these ancestors have lived in the past. Let us say that the underlying motion $\xi$ has the *distinct path property* if the law of $(\xi_s)_{s \in [0,t]}$ (considered as a $\mathcal{D}_E[0,t]$-valued random variable) is atomless for every $t > 0$ and for every initial state $\xi_0 = x \in E$. This is called Property S in [2], Definition 12.2.2.6, and occurs as formula (3.18) in [3]. In this case, the idea is that different ancestors follow a.s. different paths, and therefore it should be possible to recover the genealogy from the paths. As an immediate consequence of Lemma 17, we have the following lemma. (An analogue of this result in a spatially homogeneous setting, but for more general branching mechanisms, can be found in [3], Proposition 4.1.8(b).)

LEMMA 20 (Atomless historical superprocesses). *If $\xi$ has the distinct path property, then $\hat{\mathcal{X}}_t$ is atomless a.s. $\forall t > 0$.*

The following characterization of historical superprocesses will be convenient more than once.

LEMMA 21 (Finite-dimensional projections). *Let $\mathcal{X}$ be a $(\xi, \alpha, \beta)$-superprocess with log-Laplace semigroup $\mathcal{U} = \mathcal{U}(\xi, \alpha, \beta)$ and let $\hat{\mathcal{X}}$ be the associated historical $(\xi, \alpha, \beta)$-superprocess. Then, for all $n \geq 0$, $0 = t_0 < t_1 < \cdots < t_{n+1}$ and $f \in B_+(E^{n+2})$,*

$$\begin{aligned}E^{t_n, \hat{\mu}} &\left[ \exp\left( -\int_{\mathcal{D}_E[0,t_{n+1}]} \hat{\mathcal{X}}_{t_{n+1}}(dw) f(w_{t_0}, \ldots, w_{t_{n+1}}) \right) \right] \\ &= \exp\left( -\int_{\mathcal{D}_E[0,t_n]} \hat{\mu}(dw) \mathcal{U}_{t_{n+1}-t_n} f(w_{t_0}, \ldots, w_{t_n}, \cdot)(w_{t_n}) \right). \end{aligned} \tag{2.28}$$

Conversely, any Markov process $\hat{\mathcal{X}}$ with time-dependent state space $\mathcal{M}(\mathcal{D}_E[0,t])$ and continuous sample paths, satisfying (2.28), is the historical $(\xi,\alpha,\beta)$-superprocess.

PROOF. The fact that $\hat{\mathcal{X}}$ satisfies (2.28) can be found in [3], Theorem 2.2.5(b) or [2], Theorem 12.3.4. Conversely, if a Markov process $\hat{\mathcal{X}}$ satisfies (2.28), then, for all $0 \leq k \leq n$,

$$
(2.29) \quad \begin{aligned}
&E^{t_k,\hat{\mu}}\bigg[\exp\bigg(-\int_{\mathcal{D}_E[0,t_{n+1}]} \hat{\mathcal{X}}_{t_{n+1}}(dw) f(w_{t_0},\ldots,w_{t_{n+1}})\bigg)\bigg] \\
&= \exp\bigg(-\int_{\mathcal{D}_E[0,t_k]} \hat{\mu}(dw) f_k(w_{t_0},\ldots,w_{t_k})\bigg),
\end{aligned}
$$

where we have inductively defined functions $f_l \in B_+(E^{l+1})$ by

$$
(2.30) \quad \begin{aligned}
f_{n+1}(x_0,\ldots,x_{n+1}) &:= f(x_0,\ldots,x_{n+1}), \\
f_l(x_0,\ldots,x_l) &:= \mathcal{U}_{t_{l+1}-t_l} f_{l+1}(x_0,\ldots,x_l,\cdot)(x_l), \qquad k \leq l \leq n.
\end{aligned}
$$

The expectations in (2.29) clearly determine the transition probabilities of $\hat{\mathcal{X}}$ uniquely. □

Note that formula (2.29) says that if $\hat{\mathcal{U}}$ denotes the (time-inhomogeneous) log-Laplace semigroup of $\hat{\mathcal{X}}$ and $F(w) := f(w_{t_0},\ldots,w_{t_{n+1}})$, then

$$(2.31) \qquad \hat{\mathcal{U}}_{t_k,t_{n+1}} F(w) = f_k(w_{t_0},\ldots,w_{t_k}).$$

LEMMA 22 (Mean of historical superprocess). *Let $\hat{\mathcal{X}}$ be the historical $(\xi,\alpha,\beta)$-superprocess. Then, for any $\mu \in \mathcal{M}_1(E)$ and $m \geq 0$,*

$$(2.32) \quad E^{m\mu}[\hat{\mathcal{X}}_t](dw) = m \exp\bigg(\int_0^t \beta(w_s)\,ds\bigg) P^{\mu}[(\xi_s)_{s \in [0,t]} \in dw], \qquad t \geq 0.$$

*In particular, if $\alpha = 0$, then $\hat{\mathcal{X}}_t$ is deterministic and given by the right-hand side of (2.32).*

PROOF. By Lemma 15, the mean of a superprocess does not depend on the activity. Therefore, it suffices to prove that the historical $(\xi,0,\beta)$-superprocess is deterministic and given by the right-hand side of (2.32). Define $\hat{\mathcal{X}}_t(dw)$, $t \geq 0$, by the right-hand side of (2.32). Let $\mathcal{U} = \mathcal{U}(\xi,0,\beta)$ denote the log-Laplace semigroup of the (nonhistorical) $(\xi,0,\beta)$-superprocess. Since $\alpha = 0$, $\mathcal{U}$ coincides with the linear semigroup $\mathcal{V}$ in formula (2.12). It follows that, for $n \geq 0$, $0 = t_0 < t_1 < \cdots < t_{n+1}$ and $f \in B_+(E^{n+2})$,

$$\int_{\mathcal{D}_E[0,t_{n+1}]} \hat{\mathcal{X}}_{t_{n+1}}(dw) f(w_{t_0},\ldots,w_{t_{n+1}})$$

$$= \int_{\mathcal{D}_E[0,t_{n+1}]} m \exp\left(\int_0^{t_{n+1}} \beta(w_s)\, ds\right) f(w_{t_0},\ldots,w_{t_{n+1}})$$
$$\times P^\mu[(\xi_s)_{s\in[0,t]} \in dw]$$

(2.33)
$$= mE^\mu\left[\exp\left(\int_0^{t_{n+1}} \beta(\xi_s)\, ds\right) f(\xi_{t_0},\ldots,\xi_{t_{n+1}})\right]$$

$$= mE^\mu\left[\exp\left(\int_0^{t_n} \beta(\xi_s)\, ds\right)\right.$$
$$\left.\times E\left[\exp\left(\int_{t_n}^{t_{n+1}} \beta(\xi_s)\, ds\right) f(\xi_{t_0},\ldots,\xi_{t_{n+1}}) \big| (\xi_s)_{s\in[0,t_n]}\right]\right]$$

$$= mE^\mu\left[\exp\left(\int_0^{t_n} \beta(\xi_s)\, ds\right) \tilde{f}(\xi_{t_0},\ldots,\xi_{t_n})\right]$$

$$= \int_{\mathcal{D}_E[0,t_n]} \hat{\mathcal{X}}_{t_n}(dw)\tilde{f}(w_{t_0},\ldots,w_{t_n}),$$

where

(2.34) $$\tilde{f}(x_0,\ldots,x_n) := \mathcal{U}_{t_{n+1}-t_n} f(x_0,\ldots,x_n,\cdot)(x_n).$$

Thus, $\hat{\mathcal{X}}$ satisfies (2.28). Since $\hat{\mathcal{X}}$ is a Markov process with continuous sample paths, it follows from Lemma 21 that $(\hat{\mathcal{X}}_t)_{t\geq 0}$ is the historical $(\xi,0,\beta)$-superprocess started at time 0 in $m\mu$. $\square$

Although the next result may appear obvious, be aware of the fact that since the functions involved are not continuous, parts (b) and (c) are not trivial consequences of part (a). We will need (c) in the proof of Lemma 13.

LEMMA 23 (Preserved past property). *Let $\hat{\mathcal{X}}$ be the historical $(\xi,\alpha,\beta)$-superprocess started at time $s \geq 0$ in $\hat{\mu} \in \mathcal{M}(\mathcal{D}_E[0,s])$.*

(a) *If $F \subset \mathcal{D}_E[0,s]$ is measurable, then*

(2.35) $$P^{s,\hat{\mu}}[\hat{\mathcal{X}}_t \circ \pi_{[0,s]}^{-1} \in \mathcal{M}(F)] = 1 \qquad \forall t \geq s, \hat{\mu} \in \mathcal{M}(F).$$

(b) *If $F \subset \mathcal{D}_E[0,s]$ is a $G_\delta$-set, then*

(2.36) $$P^{s,\hat{\mu}}[\hat{\mathcal{X}}_t \circ \pi_{[0,s]}^{-1} \in \mathcal{M}(F)\ \forall t \geq s] = 1, \qquad \hat{\mu} \in \mathcal{M}(F).$$

(c) *If $F, F^c \subset \mathcal{D}_E[0,s]$ are $G_\delta$-sets, then*

(2.37) $$P^{s,\hat{\mu}}[\mathbb{1}_{\{\hat{\mathcal{X}}_{t'} \circ \pi_{[0,s]}^{-1}(F) > 0\}} \leq \mathbb{1}_{\{\hat{\mathcal{X}}_t \circ \pi_{[0,s]}^{-1}(F) > 0\}}\ \forall t' \geq t \geq s] = 1.$$

PROOF. Recall the definition of $\tilde{E}$ in (2.26) and set $\tilde{F} := \{(t, w) \in \tilde{E} : t \geq s, \pi_{[0,s]}(w) \in F\}$. If $F$ is measurable, then $\tilde{F}$ is measurable. Moreover, since $\pi_{[0,s]}$ is the pointwise limit of a sequence of continuous functions (cf. [9], Proposition 3.7.1), $\tilde{F}$ is a $G_\delta$-set when $F$ is a $G_\delta$-set. The path process $\hat{\xi}$ satisfies

$$(2.38) \quad P^{s',w}[(t, \hat{\xi}_t) \in \tilde{F} \ \forall t \geq s', (t, \hat{\xi}_{t-}) \in \tilde{F} \ \forall t > s'] = 1, \qquad (s', w) \in \tilde{F}.$$

Therefore (a) follows from Lemma 19(a) and (b) follows from Lemma 19(b). To prove (c), use the branching property (Lemma 14) to write

$$(2.39) \qquad \hat{\mathcal{X}}_t^{s,\hat{\mu}} = \hat{\mathcal{X}}_t^{s, \mathbb{1}_F \hat{\mu}} + \hat{\mathcal{X}}_t^{s, \mathbb{1}_{F^c} \hat{\mu}} \qquad \forall t \geq s \text{ a.s.}$$

Then, applying (b) to $F$ and $F^c$,

$$(2.40) \quad \begin{aligned} \hat{\mathcal{X}}_t^{s,\hat{\mu}} \circ \pi_{[0,s]}^{-1}(F) &= \hat{\mathcal{X}}_t^{s, \mathbb{1}_F \hat{\mu}} \circ \pi_{[0,s]}^{-1}(F) + \hat{\mathcal{X}}_t^{s, \mathbb{1}_{F^c} \hat{\mu}} \circ \pi_{[0,s]}^{-1}(F) \\ &= \langle \hat{\mathcal{X}}_t^{s, \mathbb{1}_F \hat{\mu}} \circ \pi_{[0,s]}^{-1}, 1 \rangle + 0 \qquad \forall t \geq s \text{ a.s.} \end{aligned}$$

By applying the strong Markov property to the stopping time $\inf\{t \geq s : \hat{\mathcal{X}}_t^{s, \mathbb{1}_F \hat{\mu}} = 0\}$, it is not hard to see that

$$(2.41) \qquad \mathbb{1}_{\{\hat{\mathcal{X}}_{t'}^{s, \mathbb{1}_F \hat{\mu}} \circ \pi_{[0,s]}^{-1} > 0\}} \leq \mathbb{1}_{\{\hat{\mathcal{X}}_t^{s, \mathbb{1}_F \hat{\mu}} \circ \pi_{[0,s]}^{-1} > 0\}} \qquad \forall t' \geq t \geq s \text{ a.s.},$$

which proves (c). □

2.2.4. *Historical binary branching particle systems.* Historical binary branching particle systems can be introduced in much the same way as historical superprocesses. First, binary branching particle systems, the underlying motion of which is a Hunt process $\xi$ with cadlag sample paths in a Polish space $E$, are defined through their generating semigroup, which in turn is defined via the unique solution to a Cauchy integral equation of the form (2.8). If $\xi$ is such a Hunt process and $b, d \in B_+(E)$, then the historical $(\xi, b, d)$-particle system $\hat{X}$ is the (time-inhomogeneous) $(\hat{\xi}, \hat{b}, \hat{d})$-particle system, where $\hat{\xi}$ is the path process associated with $\xi$ and $\hat{b}(t, w) := b(w(t))$, $\hat{d}(t, w) := d(w(t))$. Because this is very similar to what we have already seen (but easier), we skip the details.

Many of the elementary properties of historical superprocesses have analogues for historical binary branching particle systems. For example, if the underlying motion has the distinct path property, then the historical binary branching particle system at time $t > 0$ is a.s. a simple point measure. (One way to prove this is to use Poissonization and Lemma 20.) Also the formula for the finite-dimensional projections of a historical superprocess (Lemma 21) has a straightforward analogue for particle systems.

2.3. *Compensated h-transforms and weighted superprocesses.*

2.3.1. *Preliminaries from semigroup theory.* Let $E$ be a compact metrizable space and let $\mathcal{C}(E)$ be the Banach space of continuous real functions on $E$, equipped with the supremum norm, denoted by $\|\cdot\|$. Let $S = (S_t)_{t \geq 0}$ be a semigroup of bounded linear operators on $\mathcal{C}(E)$. By definition, $S$ is *strongly continuous* if $\lim_{t \to 0} \|S_t f - f\| = 0$ for all $f \in \mathcal{C}(E)$ and $S$ is *positive* if $f \geq 0$ implies $S_t f \geq 0$, $t \geq 0$. For $\lambda \in \mathbb{R}$, let us say that $S$ is $\lambda$-*contractive* if $\|S_t f\| \leq e^{\lambda t} \|f\|$, $t \geq 0$. The following version of the Hille–Yosida theorem can easily be derived from [9], Theorem 4.2.2 and Proposition 1.1.5(b). (Setting $\tilde{S}_t := e^{-\lambda} S_t$ and $\tilde{G} := G - \lambda$, we can restrict ourselves to contraction semigroups and operators $G$ that satisfy the positive maximum principle. To see that for contraction semigroups our condition (iv) implies condition (c) from [9], Theorem 4.2.2, note that $v := \int_0^\infty u_t e^{-ct} dt$ solves $(c - G)v = f$. By [9], Proposition 1.1.5(b), our condition (iv) is also necessary.)

LEMMA 24 (Hille–Yosida theorem). *A linear operator $G$ on $\mathcal{C}(E)$ with domain $\mathcal{D}(G)$ is the generator of a strongly continuous, positive, $\lambda$-contractive semigroup $S$ on $\mathcal{C}(E)$, with $\lambda \in \mathbb{R}$, if and only if*

(2.42)
- (i) *$G$ is closed;*
- (ii) *$\mathcal{D}(G)$ is dense in $\mathcal{C}(E)$;*
- (iii) *$Gf(x) \leq \lambda f(x)$ whenever $f \in \mathcal{D}(G)$ assumes its maximum over $E$ in a point $x \in E$ with $f(x) \geq 0$;*
- (iv) *for all $f \in \mathcal{D}(G)$ there exists a continuously differentiable $u : [0, \infty) \to \mathcal{C}(E)$ such that $u_0 = f$, $u_t \in \mathcal{D}(G)$ and $\frac{\partial}{\partial t} u_t = G u_t$, $t \geq 0$.*

*The function $u$ in* (iv) *is unique and given by $S_t f = u_t$, $t \geq 0$, $f \in \mathcal{D}(G)$.*

Let $G$ be the generator of a strongly continuous, positive, $\lambda$-contractive semigroup on $\mathcal{C}(E)$ and let $\alpha \in \mathcal{C}_+(E)$, $\beta \in \mathcal{C}(E)$. By definition, a *mild* solution to the Cauchy problem (1.7) is a continuous function $u : [0, \infty) \to \mathcal{C}(E)$ that satisfies

(2.43) $$u_t = S_t f + \int_0^t S_{t-s}(\beta u_s - \alpha u_s^2) \, ds, \qquad t \geq 0,$$

[cf. (2.8)]. By definition, $u$ is a *classical* solution to (1.7) if $t \mapsto u_t$ is continuously differentiable in $\mathcal{C}(E)$, $u_t \in \mathcal{D}(G)$ for all $t \geq 0$ and (1.7) holds. Every classical solution is a mild solution. For classical solutions, we have the following comparison result.

LEMMA 25 (Sub- and supersolutions). *Fix $T > 0$ and assume that $u$ is a classical solution to (1.7) on $[0, T]$ for some $u_0 = f \in \mathcal{D}(G)$. Assume that $\tilde{u} : [0, T] \to \mathcal{C}(E)$ is continuously differentiable, $\tilde{u}_t \in \mathcal{D}(G)$ for all $t \in [0, T]$ and*

$$\frac{\partial}{\partial t} \tilde{u}_t \leq G\tilde{u}_t + \beta \tilde{u}_t - \alpha \tilde{u}_t^2, \qquad t \in [0, T],$$

(2.44)

$$\tilde{u}_0 \leq f.$$

*Then $\tilde{u}_T \leq u_T$. The same holds with all inequality signs reversed.*

PROOF. This is a standard application of the maximum principle (see, e.g. [14], Lemma 10). □

Existence of solutions to (1.7) is guaranteed by the following lemma.

LEMMA 26 (Classical and mild solutions to a semilinear Cauchy problem). *For each $f \in \mathcal{C}(E)$ there exists a unique mild solution $u$ of (1.7) up to an "explosion time" $T(f)$, with $\lim_{t \uparrow T(f)} \|u_t\| = \infty$ if $T(f)$ is finite. For each $t \geq 0$, $f \mapsto \mathcal{U}_t f := u_t$ defines a continuous map from $\{f \in \mathcal{C}(E) : T(f) < t\}$, into $\mathcal{C}(E)$. If $f \in \mathcal{D}(G)$, then the mild solution to (1.7) is a classical solution. The time $T(f)$ is infinite if $f \geq 0$, in which case also $u \geq 0$, or if $\alpha = 0$.*

PROOF. The statements about mild solutions follow from [22], Theorems 6.1.2 and 6.1.4, and the statement about classical solutions follows from [22], Theorem 6.1.5. If $f \in \mathcal{D}(G) \cap \mathcal{C}_+(E)$, then using Lemma 25 it is easy to prove that the classical solution to (1.7) satisfies $0 \leq u \leq e^{(\lambda + \|\beta\|)t} \|f\|$. Since $\mathcal{D}(G)$ is dense, $\mathcal{C}_+(E)$ is the closure of its interior and $\mathcal{U}_t$ is continuous, the same bounds hold for mild solutions. The fact that solutions do not explode in the linear case $\alpha = 0$ follows from [22], Theorem 6.1.2. □

2.3.2. *Superprocesses with Feller underlying motion.* Let $E$ be a locally compact metrizable space and let $(\xi^x)^{x \in E}$ be a Markov process in $E$ with cadlag sample paths. Then $(\xi^x)^{x \in E}$ is called a *Feller process* if the map $(t, x) \mapsto \mathcal{L}(\xi_t^x)$ from $[0, \infty) \times E$ into $\mathcal{M}(E)$ is continuous and (in case $E$ is not compact) the semigroup of $(\xi^x)^{x \in E}$ maps the space $\mathcal{C}_0(E)$ of continuous real functions vanishing at infinity into itself. A Feller process on a locally compact but not compact space $E$ can always be extended to a Feller process on the one-point compactification of $E$ by putting $\xi_t^\infty := \infty$, $t \geq 0$.

If $E$ is compact, then $(\xi^x)^{x \in E}$ is a Feller process if and only if its semigroup is strongly continuous, positive, and satisfies $S_t 1 = 1$, $t \geq 0$. Such semigroups are called Feller semigroups. Note that a Feller semigroup is contractive, that is, $\lambda$-contractive with $\lambda = 0$. To every Feller semigroup there exists a unique (in law) Feller process in $E$ with cadlag sample paths ([9], Theorem 4.2.7). A Feller process on a compact metrizable space is a Hunt process (see [24], Theorem I.9.26 and Exercise I.9.27 or [16], (9.11)).

Let $E$ be compact and metrizable, let $G$ be the generator of a Feller semigroup $(P_t)_{t\geq 0}$ on $\mathcal{C}(E)$, $\alpha \in \mathcal{C}_+(E)$, and $\beta \in \mathcal{C}(E)$. Then we have the following lemma:

LEMMA 27 (Feller property of superprocess). *Let $\mathcal{X}$ be the $(G,\alpha,\beta)$-superprocess with log-Laplace semigroup $\mathcal{U} = \mathcal{U}(G,\alpha,\beta)$. Then $\mathcal{X}$ is a Feller process. For each $f \in \mathcal{C}_+(E)$, the map $(t,x) \mapsto \mathcal{U}_t f(x)$ from $[0,\infty) \times E$ into $[0,\infty)$ is continuous.*

PROOF. Since $E$ is compact, the space $\mathcal{M}(E)$ is locally compact. By [22], Theorem 6.1.4, $(t,x) \mapsto \mathcal{U}_t f(x)$ is jointly continuous in $t$ and $x$ whenever $f \in \mathcal{C}_+(E)$. Therefore, and by (1.8),

$$(2.45) \quad E^{\mu_n}[e^{-\langle \mathcal{X}_{t_n}, f \rangle}] \to E^{\mu}[e^{-\langle \mathcal{X}_t, f \rangle}] \qquad \text{as } \mu_n \Rightarrow \mu, t_n \to t, f \in \mathcal{C}_+(E).$$

If $f \in \mathcal{C}_+(E)$ satisfies $f > 0$, then the function $\mu \mapsto e^{-\langle \mu, f \rangle}$ is continuous on $\mathcal{M}(E)$ and vanishes at infinity, and by the Stone–Weierstrass theorem, the linear span of all such functions is dense in $\mathcal{C}_0(\mathcal{M}(E))$. Thus, (2.45) implies that $\mathcal{L}^{\mu_n}(\mathcal{X}_{t_n}) \Rightarrow \mathcal{L}^{\mu}(\mathcal{X}_t)$ whenever $\mu_n \Rightarrow \mu$, $t_n \to t$. It is not hard to see that the semigroup of $\mathcal{X}$ maps functions that vanish at infinity into functions that vanish at infinity; therefore, $\mathcal{X}$ is a Feller process. □

2.3.3. *Compensated h-transforms of Feller processes.* In this section we prove Lemma 3. We start with two simple observations.

LEMMA 28 (h-transformed semigroup). *Let $S$ be a strongly continuous, positive, $\lambda$-contractive semigroup on $\mathcal{C}(E)$ with generator $G$ and assume that $h \in \mathcal{D}(G)$ satisfies $h > 0$. Then*

$$(2.46) \qquad \tilde{S}_t f := \frac{1}{h} S_t(hf), \qquad f \in \mathcal{C}(E), t \geq 0,$$

*defines a strongly continuous, positive, $\tilde{\lambda}$-contractive semigroup on $\mathcal{C}(E)$, with $\tilde{\lambda} := \|\frac{Gh}{h}\|$ and generator*

$$(2.47) \quad \tilde{G}f := \frac{1}{h} G(hf) \qquad \text{with } \mathcal{D}(\tilde{G}) := \{f \in \mathcal{C}(E) : hf \in \mathcal{D}(G)\}.$$

PROOF. Since $h$ is bounded away from zero and $S$ is strongly continuous, it is easy to see that also $\tilde{S}$ is strongly continuous. Moreover, $t^{-1}(\tilde{S}_t f - f)$ converges in $\mathcal{C}(E)$ if and only if $hf \in \mathcal{D}(G)$, and the limit is given by $\tilde{G}f$. Obviously, $\tilde{S}$ is positive. Since $\frac{\partial}{\partial t} h e^{\tilde{\lambda} t} = \|\frac{Gh}{h}\| h e^{\tilde{\lambda} t} \geq Gh e^{\tilde{\lambda} t}$, Lemma 25 shows that $S_t h \leq h e^{\tilde{\lambda} t}$. Since $fh \leq \|f\|h$ and $S$ is positive, it follows that $\frac{1}{h} S_t(hf) \leq \frac{1}{h} S_t(\|f\|h) \leq \|f\| e^{\tilde{\lambda} t}$. Similarly $-\|f\| e^{\tilde{\lambda} t} \leq \tilde{S}_t f$ and, therefore, $\tilde{S}$ is $\tilde{\lambda}$-contractive. $\square$

LEMMA 29 (Linear perturbation). *Let $G$ be the generator of a strongly continuous, positive, $\lambda$-contractive semigroup on $\mathcal{C}(E)$ and assume that $g \in \mathcal{C}(E)$. Then*

$$(2.48) \qquad \tilde{G} := G + g \qquad \text{with } \mathcal{D}(\tilde{G}) := \mathcal{D}(G)$$

*is the generator of a strongly continuous, positive, $\tilde{\lambda}$-contractive semigroup on $\mathcal{C}(E)$ with $\tilde{\lambda} := \lambda + \|g\|$.*

PROOF. The operator $\tilde{G}$ satisfies conditions (i)–(iv) from Lemma 24, where condition (iv) follows from Lemma 26. $\square$

PROOF OF LEMMA 3. It follows from the previous two lemmas that $G^h$ is the generator of a strongly continuous, positive, $\lambda$-contractive semigroup on $\mathcal{C}(E)$ (for some $\lambda$). Obviously $1 \in \mathcal{D}(G^h)$ and $G^h 1 = 0$, and therefore $G^h$ generates a Feller semigroup.

To see that the law of the corresponding Feller process $\xi^h$ is given by (1.17), we proceed as follows. By [9], Lemma 4.3.2, the process

$$(2.49) \qquad M_t := \frac{h(\xi_t)}{h(x)} \exp\left( -\int_0^t \frac{Gh(\xi_s)}{h(\xi_s)} ds \right), \qquad t \geq 0,$$

is a martingale with respect to the filtration $(\mathcal{F}_t)_{t \geq 0}$ generated by $\xi$; therefore, $\tilde{P}^x(A) := E^x[M_t \mathbb{1}_A]$, $A \in \mathcal{F}_t$, defines a legitimate change of measure. Put $P_t^h f(x) := \tilde{E}^x[f(\xi_t)]$, $x \in E, f \in \mathcal{C}(E)$. We need to show that under the changed measure, $\xi$ is a Feller process with semigroup $P^h$ and that $G^h$ is the generator of $P^h$. By the Markov property of $P^x$, for $0 \leq s \leq t$,

$$(2.50) \qquad E^x\left[ f(\xi_t) \frac{M_t}{M_s} \Big| \mathcal{F}_s \right] = E^x\left[ f(\xi_t) \frac{h(\xi_t)}{h(\xi_s)} \exp\left( -\int_s^t \frac{Gh(\xi_u)}{h(\xi_u)} du \right) \Big| \mathcal{F}_s \right]$$
$$= E^{\xi_s}[f(\xi_{t-s}) M_s] = P_{t-s}^h f(\xi_s).$$

Therefore, for any $A \in \mathcal{F}_s$,

$$(2.51) \qquad \tilde{E}^x[f(\xi_t) \mathbb{1}_A] = E^x[f(\xi_t) M_t \mathbb{1}_A] = E^x\left[ E^x\left[ f(\xi_t) \frac{M_t}{M_s} \Big| \mathcal{F}_s \right] M_s \mathbb{1}_A \right]$$
$$= E^x[P_{t-s}^h f(\xi_s) M_s \mathbb{1}_A] = \tilde{E}^x[P_{t-s}^h f(\xi_s) \mathbb{1}_A],$$

which shows that $\tilde{E}^x[f(\xi_t)|\mathcal{F}_s] = P_{t-s}^h f(\xi_s)$. It is not hard to see that $P_t^h f(x)$ is jointly continuous in $t$ and $x$, and therefore $P^h$ is a Feller semigroup.

Finally, if $fh \in \mathcal{D}(G)$, then

$$h(x) \lim_{t \to 0} t^{-1}(P_t^h f - f)(x)$$

(2.52)
$$= \lim_{t \to 0} t^{-1} \left( E^x \left[ f(\xi_t) h(\xi_t) \exp\left( -\int_0^t \frac{Gh(\xi_s)}{h(\xi_s)} ds \right) \right] - h(x)f(x) \right)$$

$$= G(fh)(x) - f(x)Gh(x)$$

uniformly in $x \in E$, which shows that $G^h$ is the generator of $P^h$. □

An alternative proof of formula (1.17), using historical superprocesses, is given at the end of the next section.

2.3.4. *Weighted superprocesses.*

PROOF OF LEMMA 5. Write $\mathcal{U} := \mathcal{U}(G, \alpha, \beta)$ and $\mathcal{U}^h := \mathcal{U}(G^h, h\alpha, \beta + \frac{Gh}{h})$. By Lemma 26, for every $f \in \mathcal{D}(G^h) \cap \mathcal{C}_+(E)$, the function $t \mapsto u_t := \mathcal{U}_t(hf)$ is a classical solution to the Cauchy problem

(2.53)
$$\frac{\partial}{\partial t} u_t = Gu_t + \beta u_t - \alpha u_t^2, \qquad t \geq 0,$$
$$u_0 = hf.$$

A little calculation shows that $t \mapsto u_t^h := \frac{1}{h} u_t$ is a classical solution to the Cauchy problem

(2.54)
$$\frac{\partial}{\partial t} u_t^h = G^h u_t^h + \left( \beta + \frac{Gh}{h} \right) u_t^h - h\alpha(u_t^h)^2, \qquad t \geq 0,$$
$$u_0^h = f.$$

Therefore, $\mathcal{U}_t^h f = \frac{1}{h} \mathcal{U}_t(hf)$ for all $f \in \mathcal{D}(G^h) \cap \mathcal{C}_+(E)$. Since $\mathcal{D}(G^h)$ is dense in $\mathcal{C}(E)$, $\mathcal{C}_+(E)$ is the closure of its interior and $\mathcal{U}, \mathcal{U}^h$ are continuous, it follows that

(2.55)
$$\mathcal{U}_t^h f = \frac{1}{h} \mathcal{U}_t(hf), \qquad t \geq 0, f \in \mathcal{C}_+(E).$$

It is clear that the process $\hat{\mathcal{X}}^h$ defined in (1.18) is a Markov process with continuous sample paths. To see that $\hat{\mathcal{X}}^h$ is the historical $(G^h, h\alpha, \beta + \frac{Gh}{h})$-superprocess, by Lemma 21, it suffices to check that $\hat{\mathcal{X}}^h$ satisfies (2.28) for the log-Laplace semigroup $\mathcal{U}^h$. This is easily done, since we have

$$E\left[ \exp\left( -\int_{\mathcal{D}_E[0, t_{n+1}]} \hat{\mathcal{X}}_{t_{n+1}}^h(dw) f(w_{t_0}, \ldots, w_{t_{n+1}}) \right) \Big| \hat{\mathcal{X}}_{t_n}^h = \mu \right]$$

$$= E\left[ \exp\left( -\int_{\mathcal{D}_E[0, t_{n+1}]} h(w_{t_{n+1}}) \hat{\mathcal{X}}_{t_{n+1}}(dw) \right. \right.$$

$$\times f(w_{t_0},\ldots,w_{t_{n+1}}))\Big|(h\circ\pi_{t_n})\hat{\mathcal{X}}_{t_n}=\mu\bigg]$$

$$=E\bigg[\exp\bigg(-\int_{\mathcal{D}_E[0,t_{n+1}]}\hat{\mathcal{X}}_{t_{n+1}}(dw)h(w_{t_{n+1}})$$

(2.56)
$$\times f(w_{t_0},\ldots,w_{t_{n+1}}))\Big|\hat{\mathcal{X}}_{t_n}=(h\circ\pi_{t_n})^{-1}\mu\bigg]$$

$$=\exp\bigg(-\int_{\mathcal{D}_E[0,t_n]}h(w_{t_n})^{-1}\mu(dw)$$

$$\times\mathcal{U}_{t_{n+1}-t_n}\{h(\cdot)f(w_{t_0},\ldots,w_{t_n},\cdot)\}(w_{t_n})\bigg)$$

$$=\exp\bigg(-\int_{\mathcal{D}_E[0,t_n]}\mu(dw)\mathcal{U}^h_{t_{n+1}-t_n}f(w_{t_0},\ldots,w_{t_n},\cdot)(w_{t_n})\bigg).\quad\square$$

ALTERNATIVE PROOF OF (1.17). Let $\hat{\mathcal{X}}$ be the (deterministic) historical $(G,0,0)$-superprocess started in $\hat{\mathcal{X}}_0=\delta_x$ and set

(2.57) $$\hat{\mathcal{X}}^h_t(dw):=h(w_t)\hat{\mathcal{X}}_t(dw),\qquad t\geq 0.$$

By Lemma 5, $\hat{\mathcal{X}}^h$ is the historical $(G^h,0,\frac{Gh}{h})$-superprocess started in $\hat{\mathcal{X}}_0=h(x)\delta_x$ and, therefore, by Lemma 22,

(2.58)
(i) $\hat{\mathcal{X}}_t(dw)=P^x[(\xi_s)_{s\in[0,t]}\in dw],$

(ii) $\hat{\mathcal{X}}^h_t(dw)=h(x)\exp\bigg(\int_0^t\frac{Gh}{h}(w_s)\,ds\bigg)P^x[(\xi^h_s)_{s\in[0,t]}\in dw].$

Combining (2.57) and (2.58), we arrive at (1.17). $\square$

## 3. Proof of the main results.

### 3.1. *The infinitesimal survival probability.*

#### 3.1.1. *Extinction versus unbounded growth.*

LEMMA 30 (Eventual extinction). *We have $\mathcal{U}_t\infty\downarrow p$ as $t\uparrow\infty$. Moreover,*

(3.1) $$P^\mu[\mathcal{X}_t=0\ \text{eventually}]=\begin{cases}e^{-\langle\mu,p\rangle},&\text{if }\langle\mu,\mathcal{U}_t\infty\rangle<\infty\text{ for some }t>0,\\0,&\text{otherwise.}\end{cases}$$

*If $\sup_{x\in E}\mathcal{U}_t\infty(x)<\infty$ for some $t>0$, then $\mathcal{U}_t p=p$ for all $t\geq 0$.*

If $\langle\mu,\mathcal{U}_t\infty\rangle=\infty$ for all $t\geq 0$, then possibly $e^{-\langle\mu,p\rangle}=0$, but this need not always be the case; see Example 34.

PROOF OF LEMMA 30. Since the zero measure is an absorbing state, $\mathbb{1}_{\{\mathcal{X}_t=0\}} = \mathbb{1}_{\{\mathcal{X}_r=0 \forall r\geq t\}}$ a.s. and, therefore, $\mathbb{1}_{\{\mathcal{X}_{t_n}=0\}} \uparrow \mathbb{1}_{\{\mathcal{X}_t=0 \text{ eventually}\}}$ as $t_n \uparrow \infty$ a.s. Thus, taking the limit in (1.22), we see that $\mathcal{U}_t \infty \downarrow p$. If $\langle \mu, \mathcal{U}_t \infty \rangle < \infty$ for some $t > 0$, then $\langle \mu, \mathcal{U}_t \infty \rangle \downarrow \langle \mu, p \rangle$. Taking the limit in $P^\mu[\mathcal{X}_t = 0] = e^{-\langle \mu, \mathcal{U}_t \infty \rangle}$, we arrive at (3.1). Formula (1.8) shows that $\mathcal{U}_t$ is continuous with respect to bounded decreasing sequences. Therefore, if $\sup_{x \in E} \mathcal{U}_t \infty(x) < \infty$ for some $t > 0$, then $\mathcal{U}_t p = \mathcal{U}_t(\lim_{s \uparrow \infty} \mathcal{U}_s \infty) = \lim_{s \uparrow \infty} \mathcal{U}_{t+s} \infty = p$ for all $t \geq 0$. □

LEMMA 31 (Extinction versus unbounded growth). *If* $\sup_{x \in E} \mathcal{U}_t \infty(x) < \infty$ *for some* $t > 0$, *then*

$$(3.2) \quad P^\mu\Big[\mathcal{X}_t = 0 \text{ eventually or } \lim_{t \to \infty} \langle \mathcal{X}_t, 1 \rangle = \infty\Big] = 1, \qquad \mu \in \mathcal{M}(E),$$

*and*

$$(3.3) \qquad \lim_{t \to \infty} \mathcal{U}_t f(x) = p(x) \qquad \forall x \in E, f \in \mathcal{C}(E), f > 0.$$

PROOF. Let $(\mathcal{F}_t)_{t\geq 0}$ denote the filtration generated by $\mathcal{X}$. It follows from the right property of the process $\mathcal{X}$ [see (2.5)(i)] that $t \mapsto e^{-\langle \mathcal{X}_t, p \rangle}$ is right-continuous. By Lemma 30 and convergence of bounded right-continuous martingales,

$$(3.4) \quad e^{-\langle \mathcal{X}_t, p \rangle} = P[\mathcal{X}_s = 0 \text{ eventually}|\mathcal{F}_t] \xrightarrow[t \to \infty]{} \mathbb{1}_{\{\mathcal{X}_{s=0} \text{ eventually}\}} \qquad \text{a.s.}$$

It follows that $\langle \mathcal{X}_t, p \rangle \to \infty$ a.s. on $\{\mathcal{X}_s = 0 \text{ eventually}\}^c$. Since $\|p\| < \infty$, the same conclusion holds for $\langle \mathcal{X}_t, 1 \rangle$. □

3.1.2. *Continuity of the infinitesimal survival probability.* Even though the underlying motion has the Feller property and $\alpha, \beta$ are continuous functions, $p$ need not be continuous in general, as is illustrated by the following examples, which we give without proof.

EXAMPLE 32 (Discontinuous infinitesimal survival probability). *Let $\xi$ be the deterministic Feller process in $[-1, 1]$ given by the differential equation*

$$(3.5) \qquad \frac{\partial}{\partial t} \xi_t = 1 - (\xi_t)^2, \qquad t \geq 0.$$

*Let $\mathcal{X}$ be the superprocess in $[-1, 1]$ with underlying motion $\xi$, activity $\alpha(x) := 1$ and growth parameter $\beta(x) := -x$. Then*

$$(3.6) \qquad -\log P^{\delta_x}[\mathcal{X}_t = 0 \text{ eventually}] = \begin{cases} 1, & \text{if } x = -1, \\ 0, & \text{if } x \in (-1, 1]. \end{cases}$$

Let $\mathcal{Y}$ be the superprocess in $[-1,1]$ with underlying motion $\xi$, activity $\alpha(x) := x \vee 0$ and growth parameter $\beta(x) := x \vee 0$. Then

$$\text{(3.7)} \qquad -\log P^{\delta_x}[\mathcal{Y}_t = 0 \ eventually] = \begin{cases} \infty, & if \ x = -1, \\ 1, & if \ x \in (-1, 1]. \end{cases}$$

Nevertheless, we have the following lemma.

LEMMA 33 (Continuity of the infinitesimal survival probability). *If*

$$\sup_{x \in E} \mathcal{U}_t \infty(x) < \infty$$

*for some $t > 0$ and $\inf_{x \in E} p(x) > 0$, then $p$ is continuous.*

PROOF. Our strategy is to prove that the event that $\mathcal{X}$ becomes extinct depends in a continuous way on the path of $\mathcal{X}$ and, therefore, by the Feller property, on the initial condition. To do this, we show that by observing $\mathcal{X}$ for a finite time, we can be almost certain whether $\mathcal{X}$ becomes extinct.

Set

$$\text{(3.8)} \qquad \underline{p} := \inf_{x \in E} p(x) \quad \text{and} \quad \overline{p} := \sup_{x \in E} p(x).$$

Note that by (3.1),

$$\text{(3.9)} \qquad e^{-\langle \mu, 1 \rangle \overline{p}} \leq P^{\mu}[\mathcal{X}_t = 0 \text{ eventually}] \leq e^{-\langle \mu, 1 \rangle \underline{p}}, \qquad \mu \in \mathcal{M}(E).$$

Fix $x_0 \in E$. We will show that $p$ is continuous at $x_0$. Let $0 < c < C < \infty$ and $\varepsilon', \varepsilon'' > 0$ be arbitrary. Choose continuous functions $f_0, f_1, f_\infty$ from $[0, \infty)$ into $[0,1]$, summing up to 1, such that $\mathbb{1}_{[0,c/2]} \leq f_0 \leq \mathbb{1}_{[0,c]}$, $\mathbb{1}_{[c,C]} \leq f_1 \leq \mathbb{1}_{[c/2,2C]}$, and $\mathbb{1}_{[2C,\infty)} \leq f_\infty \leq \mathbb{1}_{[C,\infty)}$. By Lemma 31, there exists a $T > 0$ such that

$$\text{(3.10)} \qquad E^{\delta_{x_0}}[f_1(\langle \mathcal{X}_T, 1 \rangle)] \leq \varepsilon'.$$

Let $d$ be a metric that generates the topology on $E$. By Lemma 27, we can choose $\delta > 0$ such that for all $x \in E$ with $d(x, x_0) \leq \delta$,

$$\text{(3.11)} \quad |E^{\delta_{x_0}}[f_r(\langle \mathcal{X}_T, 1 \rangle)] - E^{\delta_x}[f_r(\langle \mathcal{X}_T, 1 \rangle)]| \leq \varepsilon'', \qquad d(x, x_0) \leq \delta, r = 0, 1.$$

Write

$$\text{(3.12)} \qquad \begin{aligned} P^{\delta_x}[\mathcal{X}_t &= 0 \text{ eventually}] \\ &= E^{\delta_x}\left[ \sum_{r=0,1,\infty} f_r(\langle \mathcal{X}_T, 1 \rangle) \mathbb{1}_{\{\mathcal{X}_{t=0} \text{ eventually}\}} \right] \\ &= \sum_{r=0,1,\infty} E^{\delta_x}[f_r(\langle \mathcal{X}_T, 1 \rangle) P^{\mathcal{X}_T}[\mathcal{X}_t = 0 \text{ eventually}]]. \end{aligned}$$

Using (3.12) to get lower and upper estimates on $P^{\delta_x}[\mathcal{X}_t = 0 \text{ eventually}]$, and applying (3.9), we find that

$$
\begin{aligned}
& E^{\delta_x}[f_0(\langle \mathcal{X}_T, 1\rangle)] - (1 - e^{-c\overline{p}}) \\
& \quad \leq E^{\delta_x}[f_0(\langle \mathcal{X}_T, 1\rangle)] e^{-c\overline{p}} \\
& \quad \leq P^{\delta_x}[\mathcal{X}_t = 0 \text{ eventually}] \\
& \quad \leq E^{\delta_x}[f_0(\langle \mathcal{X}_T, 1\rangle)] + E^{\delta_x}[f_1(\langle \mathcal{X}_T, 1\rangle)] + E^{\delta_x}[f_\infty(\langle \mathcal{X}_T, 1\rangle)] e^{-C\underline{p}} \\
& \quad \leq E^{\delta_x}[f_0(\langle \mathcal{X}_T, 1\rangle)] + (\varepsilon' + \varepsilon'') + e^{-C\underline{p}}, \qquad d(x, x_0) \leq \delta.
\end{aligned}
$$
(3.13)

Therefore, for all $x \in E$ with $d(x, x_0) \leq \delta$,

$$
\begin{aligned}
|P^{\delta_{x_0}}[\mathcal{X}_t & = 0 \text{ eventually}] - P^{\delta_x}[\mathcal{X}_t = 0 \text{ eventually}]| \\
& \leq |P^{\delta_{x_0}}[\mathcal{X}_t = 0 \text{ eventually}] - E^{\delta_{x_0}}[f_0(\langle \mathcal{X}_T, 1\rangle)]| \\
& \quad + |E^{\delta_{x_0}}[f_0(\langle \mathcal{X}_T, 1\rangle)] - E^{\delta_x}[f_0(\langle \mathcal{X}_T, 1\rangle)]| \\
& \quad + |E^{\delta_x}[f_0(\langle \mathcal{X}_T, 1\rangle)] - P^{\delta_x}[\mathcal{X}_t = 0 \text{ eventually}]| \\
& \leq ((1 - e^{-c\overline{p}}) + (\varepsilon' + \varepsilon'' + e^{-C\underline{p}})) \\
& \quad + \varepsilon'' + ((1 - e^{-c\overline{p}}) + (\varepsilon' + \varepsilon'' + e^{-C\underline{p}})).
\end{aligned}
$$
(3.14)

Since $0 < c < C < \infty$ and $\varepsilon', \varepsilon'' > 0$ are arbitrary, the last line of (3.14) can be made arbitrarily small. Thus, we have shown that for each $\varepsilon > 0$ there exists a $\delta > 0$ such that

(3.15) $\qquad |e^{-p(x_0)} - e^{-p(x)}| \leq \varepsilon \qquad \forall x \in E \text{ with } d(x, x_0) \leq \delta.$

This shows that $p$ is continuous at $x_0$. $\square$

### 3.1.3. *Properties of the infinitesimal survival probability.*

PROOF OF PROPOSITION 7. Parts (a) and (b) follow from Lemmas 30 and 31. To prove part (c), note that if $f \in \mathcal{C}_+(E)$ satisfies $\mathcal{U}_t f = f$ for all $t \geq 0$, then $u_t := f$, $t \geq 0$, is a mild solution to (1.7), that is,

(3.16) $\qquad f = P_t f + \int_0^t P_s(\beta f - \alpha f^2)\, ds, \qquad t \geq 0.$

Thus,

(3.17) $\quad \lim_{t \to 0} t^{-1}(P_t f - f) = -\lim_{t \to 0} t^{-1} \int_0^t P_s(\beta f - \alpha f^2)\, ds = -\beta f + \alpha f^2,$

which proves that $f \in \mathcal{D}(G)$ and that (1.23) holds. Conversely, if $f \in \mathcal{D}(G) \cap \mathcal{C}_+(E)$ solves (1.23), then $u_t := f$ is a classical solution to (1.7) and, therefore, $\mathcal{U}_t f = f$ for all $t \geq 0$.

To prove (d), note that if $\inf_{x \in E} p(x) > 0$, then $p$ is continuous by Lemma 33 and, therefore, $p$ solves (1.23) by parts (b) and (c). Moreover, part (a) shows that in this case there exists only one positive fixed point of $\mathcal{U}$. $\square$

3.1.4. *Nonuniform convergence of $\mathcal{U}_t\infty$.* Lemma 11 shows that the assumption that $\sup_{x\in E}\mathcal{U}_t\infty(x) < \infty$ for some $t > 0$ cannot be dropped from Theorems 8 and 9. However, the reader may wonder if this condition is not implied by the simpler-looking condition $\sup_{x\in E} p(x) < \infty$. To show that this is not the case, we include the following example.

EXAMPLE 34 (Nonuniform convergence of $\mathcal{U}_t\infty$). *There exists a generator $G$ of a Feller process in a compact metrizable space $E$ and $\alpha \in \mathcal{C}_+(E)$ such that $\mathcal{U} = \mathcal{U}(G, \alpha, 0)$ satisfies*

(3.18)
$$\begin{aligned}&\text{(i) } \mathcal{U}_t\infty(x) < \infty &&\forall x \in E, t > 0,\\&\text{(ii) } \mathcal{U}_t\infty \downarrow 0, &&\text{as } t\uparrow\infty,\\&\text{(iii) } \sup_{x\in E}\mathcal{U}_t\infty(x) = \infty &&\forall t \geq 0.\end{aligned}$$

PROOF. Take $E := [0,1]^2$. Define a Feller process $\xi = (\xi^x)^{x\in E}$ in $E$ by

(3.19)
$$\xi_t^{(x,y)} := (x, ye^{-t}), \qquad (x,y) \in [0,1] \times [0,1),$$
$$\xi_t^{(x,1)} := \begin{cases} (x,1), & t \leq \tau_x, \\ (x, e^{-(t-\tau_x)}), & t > \tau_x, \end{cases} \qquad x \in [0,1],$$

where $\tau_x$, $x \in (0,1]$ is an exponentially distributed random variable with mean $x$ and $\tau_0 := 0$. It is not hard to see that $\xi$ is a Feller process. Let $G$ denote its generator. Choose $\alpha \in \mathcal{C}_+(E)$ such that $\alpha(0,1) = 0$ and $\alpha > 0$ elsewhere. Set

(3.20) $$\underline{\alpha(x,\cdot)} := \inf\{\alpha(x,y) : y \in [0,1]\}, \qquad x \in [0,1].$$

For fixed $x \in [0,1]$, the process $\xi$ restricted to $\{x\} \times [0,1]$ is an autonomous Feller process and $\underline{\alpha(x,\cdot)} > 0$ for $x > 0$. Therefore, using (1.30), we have

(3.21) $$\mathcal{U}_t\infty(x,y) \leq \frac{1}{\underline{\alpha(x,\cdot)}t}, \qquad t > 0, (x,y) \in (0,1] \times [0,1].$$

The superprocess $\mathcal{X}$ started in $\delta_{(0,y)}$ ($y \in [0,1]$) is concentrated on $(0, ye^{-t})$ at time $t$, if it survives. Therefore, applying (1.30) to the process $(\mathcal{X}_t)_{t\geq\varepsilon}$, we have for each $\varepsilon > 0$ that

(3.22) $$\mathcal{U}_t\infty(0,y) \leq \frac{1}{\delta(t-\varepsilon)}, \qquad t > \varepsilon, \text{ where } \delta := \inf\{\alpha(0,e^{-t}) : t \in [\varepsilon,\infty]\}.$$

This proves (3.18)(i) and (3.18)(ii). Now consider the process $(\mathcal{X}_t(\cdot\cap((0,1]\times\{1\})))_{t\geq 0}$. It is not too hard to see that this is an autonomous superprocess without (i.e., with constant) underlying motion, activity $\alpha(\cdot,1)$ and growth parameter $\beta(x) := -\frac{1}{x}$. Therefore [see (1.30)],

(3.23) $$\mathcal{U}_t(\infty\mathbb{1}_{(0,1]\times\{1\}})(x,1) = \frac{\beta(x)}{\alpha(x,1)(1-e^{-\beta(x)t})} = \frac{x^{-1}}{\alpha(x,1)(e^{t/x}-1)},$$

$t > 0, x \in (0,1]$. We can additionally choose $\alpha(x,1) := e^{-1/x^2}$, $x \in (0,1]$. Then
$$(3.24) \qquad \lim_{x \to 0} \mathcal{U}_t(\infty \mathbb{1}_{(0,1] \times \{1\}})(x,1) = \infty, \qquad t > 0.$$
It follows that $\sup_{x \in E} \mathcal{U}_t \infty(x) \geq \sup_{x \in E} \mathcal{U}_t(\infty \mathbb{1}_{(0,1] \times \{1\}})(x) = \infty$, which proves (3.18)(iii). □

### 3.2. Surviving lines of descent.

3.2.1. *Poisson point measures.* Let $E$ be a Polish space. By definition, a Poisson point measure with intensity $\mu \in \mathcal{M}(E)$ is an $\mathcal{N}(E)$-valued random variable $\mathrm{Pois}(\mu)$ with
$$(3.25) \qquad E[(1-f)^{\mathrm{Pois}(\mu)}] = e^{-\langle \mu, f \rangle}, \qquad f \in B_+(E).$$
If $\mu$ is atomless, then $\mathrm{Pois}(\mu)$ a.s. takes values in the space $\mathcal{N}^*(E) := \{\nu \in \mathcal{N}(E) : \nu(\{x\}) \leq 1 \,\forall x \in E\}$ of simple point measures on $E$. Note that $\mathcal{N}^*(E)$ is an open subset of $\mathcal{N}(E)$ and, therefore, a Polish space in the induced topology. We identify $\mathcal{N}^*(E)$ with the space of finite subsets of $E$. If $\mu \in \mathcal{M}(E)$ is atomless, then an $\mathcal{N}^*(E)$-valued random variable $\nu$ is a Poisson point measure with intensity $\mu$ if and only if (see [21], Proposition 1.4.7)
$$(3.26) \qquad P[\nu(A) = 0] = e^{-\mu(A)}, \qquad A \in \mathcal{B}(E).$$
It is not hard to see that the event $\{\mu \in \mathcal{M}(E) : \mathrm{supp}(\mu) \text{ is finite}\} \subset \mathcal{M}(E)$ is measurable and that $\mu \mapsto \mathrm{supp}(\mu)$ is a measurable map from $\{\mu \in \mathcal{M}(E) : \mathrm{supp}(\mu) \text{ is finite}\}$ into $\mathcal{N}^*(E)$.

We need a criterion to decide whether the support of a random measure is a Poisson point measure.

LEMMA 35 (Random measures with Poisson support). *Let $E$ be a Polish space, let $\mu$ be an atomless measure on $E$ and let $\mathcal{Z}$ be an $\mathcal{M}(E)$-valued random variable such that*
$$(3.27) \qquad P[\mathcal{Z}(A) = 0] = e^{-\mu(A)}, \qquad A \in \mathcal{B}(E).$$
*Then*
$$(3.28) \qquad P[\mathrm{supp}(\mathcal{Z}) \text{ is finite}] = \begin{cases} 1, & \text{if } \mu(E) < \infty, \\ 0, & \text{if } \mu(E) = \infty. \end{cases}$$
*Moreover, if $\mu(E) < \infty$, then $\mathrm{supp}(\mathcal{Z})$ is a Poisson point measure with intensity $\mu$.*

PROOF. Assume that $\mu(E) < \infty$. Choose finite measurable partitions $A^{(n)} = \{A_i^{(n)}\}_{i \in I^{(n)}}$ such that $A^{(n+1)}$ is a refinement of $A^{(n)}$ and such that intersections of the form $\bigcap A_{i_n}^{(n)}$ are empty or consist of one point. Since
$$(3.29) \qquad E[|\{i \in I^{(n)} : \mathcal{Z}(A_i^{(n)}) > 0\}|] = \sum_{i \in I^{(n)}} (1 - e^{-\mu(A_i^{(n)})}) \leq \mu(E),$$

the increasing limit of $|\{i \in I^{(n)} : \mathcal{Z}(A_i^{(n)}) > 0\}|$ is a.s. finite, that is, there are a.s. finitely many decreasing sequences of partition elements $A_{i_1}^{(1)} \supset A_{i_2}^{(2)} \supset \cdots$ such that $\mathcal{Z}(A_{i_n}^{(n)}) > 0$ for all $n$. The limit points of these sequences give the support of $\mathcal{Z}$ and by formula (3.26), $\mathrm{supp}(\mathcal{Z})$ is a Poisson point measure with intensity $\mu$.

Assume, on the other hand, that $\mu(E) = \infty$. Since $\mu$ is atomless, there exist measurable disjoint sets $(B_i)_{i \geq 0}$ such that $\mu(B_i) \geq 1$. Formula (3.27) shows that the events $\{\mathcal{Z}(B_i) > 0\}$ are independent and that

$$(3.30) \qquad \sum_{i=1}^{\infty} P[\mathcal{Z}(B_i) > 0] = \sum_{i=1}^{\infty} (1 - e^{-\mu(B_i)}) = \infty.$$

Therefore, by the Borel–Cantelli lemma $\mathcal{Z}(B_i) > 0$ for infinitely many $i$, which proves that $\mathrm{supp}(\mathcal{Z})$ is infinite a.s. □

3.2.2. *Poissonization of historical superprocesses.* The following lemma gives a historical variant of formula (1.11)(i). Moreover, it shows that the particles in $\mathrm{Pois}((\mathcal{U}_t f)\mu)$ from (1.11)(i) are, in a sense, the ancestors of the particles in $\mathrm{Pois}(f\mathcal{X}_t)$.

LEMMA 36 (Poissonization of historical superprocesses). *Let $\hat{\mathcal{X}}$ be the historical $(G, \alpha, \beta)$-superprocess started at time $s \geq 0$ in $\hat{\mu} \in \mathcal{M}(\mathcal{D}_E[0,s])$. Assume that $\hat{\mu}$ is atomless. If $\hat{\nu}$ is an $\mathcal{N}(\mathcal{D}_E[0,s+t])$-valued random variable such that, for a given $f \in B_+(E)$ and $t \geq 0$,*

$$(3.31) \quad P[\hat{\nu} \in \cdot | (\hat{\mathcal{X}}_r)_{s \leq r \leq s+t}] = P[\mathrm{Pois}((f \circ \pi_{s+t})\hat{\mathcal{X}}_{s+t}) \in \cdot | \hat{\mathcal{X}}_{s+t}] \qquad a.s.,$$

*then $\mathrm{supp}(\hat{\nu} \circ \pi_{[0,s]}^{-1})$ is a Poisson point measure with intensity $(\mathcal{U}_t f \circ \pi_s)\hat{\mu}$.*

PROOF. Since $\hat{\mu}$ is atomless, by Lemma 35, it suffices to show that for all $A \in \mathcal{B}(\mathcal{D}_E[0,s])$,

$$(3.32) \qquad P[\hat{\nu} \circ \pi_{[0,s]}^{-1}(A) = 0] = \exp\bigl(-(\mathcal{U}_t f \circ \pi_s)\hat{\mu}(A)\bigr).$$

By (3.31),

$$(3.33) \quad P[\hat{\nu} \circ \pi_{[0,s]}^{-1}(A) = 0] = E^{s,\hat{\mu}}[\exp\bigl(-(f \circ \pi_{s+t})\hat{\mathcal{X}}_{s+t} \circ \pi_{[0,s]}^{-1}(A)\bigr)].$$

By the branching property (Lemma 14) and by Lemma 23(a), we can rewrite the right-hand side of this equation as

$$
\begin{aligned}
(3.34) \quad & E^{s,\mathbb{1}_A \hat{\mu}}[\exp\bigl(-(f \circ \pi_{s+t})\hat{\mathcal{X}}_{s+t} \circ \pi_{[0,s]}^{-1}(A)\bigr)] \\
& \qquad \times E^{s,\mathbb{1}_{A^c} \hat{\mu}}[\exp\bigl(-(f \circ \pi_{s+t})\hat{\mathcal{X}}_{s+t} \circ \pi_{[0,s]}^{-1}(A)\bigr)] \\
& \quad = E^{s,\mathbb{1}_A \hat{\mu}}[\exp\bigl(-\langle (f \circ \pi_{s+t})\hat{\mathcal{X}}_{s+t}, 1\rangle\bigr)] \cdot 1.
\end{aligned}
$$

TRIMMED TREES AND EMBEDDED SYSTEMS 37

From the relation (2.27) between a historical superprocess and its associated superprocess it is obvious that

$$(3.35) \quad E^{s,\mathbb{1}_A \hat{\mu}}[\exp(-\langle \hat{\mathcal{X}}_{s+t} \circ \pi_{s+t}^{-1}, f \rangle)] = \exp(-\langle (\mathbb{1}_A \hat{\mu}) \circ \pi_s^{-1}, \mathcal{U}_t f \rangle).$$

It follows that

$$(3.36) \quad E^{s,\hat{\mu}}[\exp(-(f \circ \pi_{s+t})\hat{\mathcal{X}}_{s+t} \circ \pi_{[0,s]}^{-1}(A))] = \exp(-(\mathcal{U}_t f \circ \pi_s)\hat{\mu}(A)).$$

Combining this with (3.33), we see that (3.32) holds. □

The proof of Lemma 36 has the following corollary.

COROLLARY 37 (Surviving lines of descent). *Let $\hat{\mathcal{X}}$ be the historical $(G, \alpha, \beta)$-superprocess started at time $s \geq 0$ in $\hat{\mu} \in \mathcal{M}(\mathcal{D}_E[0,s])$. Assume that $\hat{\mu}$ is atomless. Then, for any $t > 0$,*

$$(3.37) \quad P[\mathrm{supp}(\hat{\mathcal{X}}_{s+t} \circ \pi_{[0,s]}^{-1}) \text{ is finite}] = 1 \iff \langle \hat{\mu} \circ \pi_s^{-1}, \mathcal{U}_t \infty \rangle < \infty.$$

*Moreover, if $\langle \hat{\mu} \circ \pi_s^{-1}, \mathcal{U}_t \infty \rangle < \infty$, then $\mathrm{supp}(\hat{\mathcal{X}}_{s+t} \circ \pi_{[0,s]}^{-1})$ is a Poisson point measure with intensity $(\mathcal{U}_t \infty \circ \pi_s)\hat{\mu}$.*

PROOF. Letting $f \uparrow \infty$ in (3.36) we see that

$$(3.38) \quad \begin{aligned} P^{s,\hat{\mu}}[\hat{\mathcal{X}}_{s+t} \circ \pi_{[0,s]}^{-1}(A) = 0] \\ = \exp(-(\mathcal{U}_t \infty \circ \pi_s)\hat{\mu}(A)), \quad A \in \mathcal{B}(\mathcal{D}_E[0,s]). \end{aligned}$$

Now the statements follow from Lemma 35. □

3.2.3. *Finite ancestry property.*

PROOF OF LEMMA 11. If $\sup_{x \in E} \mathcal{U}_t \infty(x) < \infty$ for some $t > 0$, then $\langle \mu, \mathcal{U}_t \infty \rangle < \infty$ for all $\mu \in \mathcal{M}(E)$. On the other hand, if $\sup_{x \in E} \mathcal{U}_t \infty(x) = \infty$ for all $t \geq 0$, then we can find $\mu \in \mathcal{M}(E)$ such that $\langle \mu, \mathcal{U}_t \infty \rangle = \infty$ for all $t \geq 0$. To see this, choose strictly positive $(\varepsilon_n)_{n \geq 0}$ such that $\sum_{n \geq 0} \varepsilon_n = 1$. Choose $t_n \uparrow \infty$ and $x_n \in E$ such that $\mathcal{U}_{t_n} \infty(x_n) \geq \varepsilon_n^{-1}$ and choose $\mu := \sum_{n \geq 0} \varepsilon_n \delta_{x_n}$. Then $\langle \mu, \mathcal{U}_{t_n} \infty \rangle \geq \sum_{m \geq n} \varepsilon_m \mathcal{U}_{t_n} \infty(x_m) \geq \sum_{m \geq n} \varepsilon_m \mathcal{U}_{t_m} \infty(x_m) = \infty$.

The log-Laplace semigroup $\mathcal{U}' = \mathcal{U}(G', \alpha', \beta')$ satisfies $\mathcal{U}'_t(f \circ \psi) = (\mathcal{U}_t f) \circ \psi$, where $\psi$ denotes the projection from $E'$ to $E$ (see Lemma 18). Therefore (i) implies that $\langle \mu \otimes \ell, \mathcal{U}'_t \infty \rangle < \infty$ for some $t > 0$, which by Corollary 37 implies (ii). On the other hand, if (i) does not hold, then there exists a $\mu \in \mathcal{M}(E)$ such that $\langle \mu \otimes \ell, \mathcal{U}'_t \infty \rangle = \infty$ for all $t \geq 0$, and in this case Corollary 37 shows that (ii) does not hold. Finally, since $\hat{\mathcal{X}}_t = \hat{\mathcal{X}}'_t \circ \psi_t^{-1}$, (ii) implies (iii). □

PROOF OF LEMMA 13. We prove the following, slightly more general result.

LEMMA 38 (Immortal lines of descent). *Let $\hat{\mathcal{X}}$ be the historical $(G, \alpha, \beta)$-superprocess started at time 0 in $\mu \in \mathcal{M}(E)$. Assume that $\sup_{x \in E} \mathcal{U}_t \infty(x) < \infty$ for all $t > q$, for some $q \geq 0$, where $\mathcal{U} = \mathcal{U}(G, \alpha, \beta)$. Then*

(3.39)
(i) $\operatorname{supp}(\hat{\mathcal{X}}_r \circ \pi_{[0,t]}^{-1})$ *is finite*

$\forall\, t, r \geq 0$ *such that* $t + q < r$ *a.s.*

(ii) $\operatorname{supp}(\hat{\mathcal{X}}_r \circ \pi_{[0,t]}^{-1}) \supset \operatorname{supp}(\hat{\mathcal{X}}_{r'} \circ \pi_{[0,t]}^{-1})$

$\forall\, t, r, r' \geq 0$ *such that* $t + q < r \leq r'$ *a.s.*

PROOF. Let us introduce the shorthand

(3.40) $$\hat{\mathcal{X}}_{t,r} := \hat{\mathcal{X}}_r \circ \pi_{[0,t]}^{-1}, \qquad 0 \leq t \leq r.$$

Let $D \subset [0, \infty)$ be countable and dense. The implication $\Leftarrow$ in (3.37) also holds if $\hat{\mu}$ is not atomless; this can be proved by extending the space $E$ as in Lemma 11. Therefore,

(3.41) $\operatorname{supp}(\hat{\mathcal{X}}_{t,r})$ is finite $\quad \forall\, t, r \in D, t + q < r$ a.s.

Let $\mathcal{O}$ be a countable basis for the topology on $\mathcal{D}_E[0,t]$. Conditioning on $\hat{\mathcal{X}}_t$ and applying Lemma 23(c), we see that

(3.42) $$\mathbb{1}_{\{\hat{\mathcal{X}}_{t,r'}(O) > 0\}} \leq \mathbb{1}_{\{\hat{\mathcal{X}}_{t,r}(O) > 0\}} \qquad \forall\, r, r' \geq 0, t \in D, O \in \mathcal{O}, t \leq r \leq r' \text{ a.s.}$$

It follows that

(3.43) $\operatorname{supp}(\hat{\mathcal{X}}_{t,r'}) \subset \operatorname{supp}(\hat{\mathcal{X}}_{t,r}) \quad \forall\, r, r' \geq 0, t \in D, t \leq r \leq r'$ a.s.

Combining this with (3.41), we see that $\operatorname{supp}(\hat{\mathcal{X}}_{t,r'}) \subset \operatorname{supp}(\hat{\mathcal{X}}_{t,r})$ and $\operatorname{supp}(\hat{\mathcal{X}}_{t,r})$ is finite $\forall\, r' \geq 0, t, r \in D, t + q < r \leq r'$ a.s., and therefore (3.41) can be sharpened to

(3.44) $\operatorname{supp}(\hat{\mathcal{X}}_{t,r'})$ is finite $\quad \forall\, r' \geq 0, t \in D, t + q < r'$ a.s.

If $\hat{\mathcal{X}}_{t,r'}$ is finitely supported for some $t, r'$, then $\operatorname{supp}(\hat{\mathcal{X}}_{t',r'}) = \pi_{[0,t']}(\operatorname{supp}(\hat{\mathcal{X}}_{t,r'}))$ for all $t' \leq t$. Thus, (3.44) can be further sharpened to

(3.45) $\operatorname{supp}(\hat{\mathcal{X}}_{t',r'})$ is finite $\quad \forall\, t', r' \geq 0, t' + q < r'$ a.s.

This proves (3.39)(i). Moreover, by (3.43) and (3.45),

(3.46) $\operatorname{supp}(\hat{\mathcal{X}}_{t',r'}) = \pi_{[0,t']}(\operatorname{supp}(\hat{\mathcal{X}}_{t,r'})) \subset \pi_{[0,t']}(\operatorname{supp}(\hat{\mathcal{X}}_{t,r})) = \operatorname{supp}(\hat{\mathcal{X}}_{t',r})$

$\forall\, t', r, r' \geq 0, t \in D, t' + q \leq t + q < r \leq r'$ a.s.,

which proves (3.39)(ii). $\square$

The proof of Lemma 13 is complete. $\square$

3.3. *Embedded trees.* Our first and crucial proposition in this section shows that it is possible to embed a collection $I$ of immortal lines of descent in certain historical superprocesses. We then identify these immortal lines of descent as a historical binary branching particle system. Finally, we generalize our results in a number of steps, until we arrive at the statements in Section 1.5.

3.3.1. *Construction of the embedded tree.* Recall the definition of the distinct path property before Lemma 20.

PROPOSITION 39 (Embedded tree). *Let $\hat{\mathcal{X}}$ be the historical $(G,\alpha,\alpha)$-superprocess started at time 0 in $\mu \in \mathcal{M}(E)$. Assume that $\mu$ is atomless and that the Feller process with generator $G$ has the distinct path property. Then $\hat{\mathcal{X}}$ may be coupled to a random set $I \subset \mathcal{D}_E[0,\infty)$ such that the random sets $I_t := \{\pi_{[0,t]}(w) : w \in I\}$ are finite for all $t \geq 0$ and satisfy*

$$(3.47) \qquad P[I_t \in \cdot | (\hat{\mathcal{X}}_s)_{0 \leq s \leq t}] = P[\mathrm{Pois}(\hat{\mathcal{X}}_t) \in \cdot | \hat{\mathcal{X}}_t] \qquad a.s. \ \forall t \geq 0.$$

*If, in addition, $\mathcal{U} = \mathcal{U}(G,\alpha,\alpha)$ satisfies $\sup_{x \in E} \mathcal{U}_t \infty(x) < \infty$ for some $t > 0$, then $p := \lim_{t \uparrow \infty} \mathcal{U}_t \infty = 1$ and $I$ may be chosen such that, moreover,*

$$(3.48) \qquad I_t = \mathrm{supp}\,(\hat{\mathcal{X}}_r \circ \pi_{[0,t]}^{-1}), \qquad r\text{-eventually } \forall t \geq 0 \ a.s.$$

PROOF. Identify, as usual, finite subsets and simple point measures. For each $T \geq 0$, let $I^{(T)}$ be a random finite subset of $\mathcal{D}_E[0,T]$ such that

$$(3.49) \qquad P[I^{(T)} \in \cdot | (\hat{\mathcal{X}}_t)_{0 \leq t \leq T}] = P[\mathrm{Pois}(\hat{\mathcal{X}}_T) \in \cdot | \hat{\mathcal{X}}_T].$$

Put

$$(3.50) \quad I_t^{(T)} := \{\pi_{[0,t]}(w) : w \in I^{(T)}\} = \mathrm{supp}\,(I^{(T)} \circ \pi_{[0,t]}^{-1}), \qquad 0 \leq t \leq T.$$

Using the fact that, by Lemma 20, $\hat{\mathcal{X}}_t$ is a.s. atomless, conditioning on $(\hat{\mathcal{X}}_s)_{0 \leq s \leq t}$, applying Lemma 36 and the fact that the function 1 is a fixed point of $\mathcal{U}(G,\alpha,\alpha)$, we find that

$$(3.51) \quad P[I_t^{(T)} \in \cdot | (\hat{\mathcal{X}}_s)_{0 \leq s \leq t}] = P[\mathrm{Pois}(\hat{\mathcal{X}}_t) \in \cdot | \hat{\mathcal{X}}_t] \qquad \text{a.s. } \forall 0 \leq t \leq T.$$

Thus, we can satisfy (3.47) up to a finite time horizon $T$. To let $T \uparrow \infty$, we need to take a projective limit. For $0 \leq S \leq T$, define a map $\psi_{S,T} : \mathcal{N}^*(\mathcal{D}_E[0,T]) \to \mathcal{N}^*(\mathcal{D}_E[0,S])$ by

$$(3.52) \qquad \psi_{S,T}(J) := \{\pi_{[0,S]}(w) : w \in J\}, \qquad J \in \mathcal{N}^*(\mathcal{D}_E[0,T]).$$

Then (3.51) shows that the random variables $((\hat{\mathcal{X}}_t)_{0 \leq t \leq T}, I^{(T)})_{T \geq 0}$ satisfy the consistency relation $\mathcal{L}((\hat{\mathcal{X}}_t)_{0 \leq t \leq S}, \psi_{S,T}(I^{(T)})) = \mathcal{L}((\hat{\mathcal{X}}_t)_{0 \leq t \leq S}, I^{(S)})$ $(0 \leq S \leq$ $T$). Note that $((\hat{\mathcal{X}}_t)_{0\leq t\leq T}, I^{(T)})$ takes values in the Polish space $\mathcal{C}_{\mathcal{M}(\mathcal{D}_E[0,\infty))}[0,T] \times \mathcal{N}^*(\mathcal{D}_E[0,T])$. Let $\mathcal{N}^{(\infty)}$ be the space of all countable subsets $I \subset \mathcal{D}_E[0,\infty)$ such that $\psi_{T,\infty}(I) := \{\pi_{[0,T]}(w) : w \in I\}$ is finite for all $T \geq 0$. Equip $\mathcal{N}^{(\infty)}$ with the $\sigma$-field generated by the mappings $\psi_{T,\infty} : \mathcal{N}^{(\infty)} \to \mathcal{N}^*(\mathcal{D}_E[0,T])$, $T \geq 0$. Taking the projective limit of the variables $((\hat{\mathcal{X}}_t)_{0\leq t\leq T}, I^{(T)})_{T\geq 0}$, we can construct a random variable $(\tilde{\mathcal{X}}, I)$ with values in $\mathcal{C}_{\mathcal{M}(\mathcal{D}_E[0,\infty))}[0,\infty) \times \mathcal{N}^{(\infty)}$ such that $((\tilde{\mathcal{X}}_t)_{0\leq t\leq T}, \psi_{T,\infty}(I))$ is equal in distribution to $((\hat{\mathcal{X}}_t)_{0\leq t\leq T}, I^{(T)})$ for all $T \geq 0$. It follows that $\tilde{\mathcal{X}}$ is the historical $(G, \alpha, \alpha)$-superprocess started at time 0 in $\mu \in \mathcal{M}(E)$ and that $I$ is a random set that satisfies (3.47).

Assume that $\sup_{x \in E} \mathcal{U}_t \infty(x) < \infty$ for some $t > 0$. We must show that we can choose $I$ such that, moreover, (3.48) holds. First note that the function 1 is a positive solution to (1.23) and, therefore, by Proposition 7(a), $p = 1$. Choose $q \geq 0$ such that $\sup_{x \in E} \mathcal{U}_t \infty(x) < \infty$ for all $t > q$. Then, by Lemma 38, the random sets $\text{supp}(\hat{\mathcal{X}}_r \circ \pi_{[0,t]}^{-1})$ are finite and nonincreasing in $r > t + q$ for all $t \geq 0$ a.s. Define random finite subsets $I_t \subset \mathcal{D}_E[0,t]$ by

$$(3.53) \qquad I_t := \bigcap_{r > t+q} \text{supp}(\hat{\mathcal{X}}_r \circ \pi_{[0,t]}^{-1}) \qquad \forall t \geq 0 \text{ a.s.}$$

Then (3.48) is fulfilled. Define $I \subset \mathcal{D}_E[0,\infty)$ by

$$(3.54) \qquad I := \{w \in \mathcal{D}_E[0,\infty) : \pi_{[0,t]}(w) \in I_t \ \forall t \geq 0\}.$$

Then

$$(3.55) \qquad I_t = \{\pi_{[0,t]}(w) : w \in I\} \qquad \forall t \geq 0 \text{ a.s.}$$

By Corollary 37,

$$(3.56) \qquad \begin{aligned} &P[\text{supp}(\hat{\mathcal{X}}_r \circ \pi_{[0,t]}^{-1}) \in \cdot | (\hat{\mathcal{X}}_s)_{0 \leq s \leq t}] \\ &\qquad = P[\text{Pois}((\mathcal{U}_{r-t}\infty \circ \pi_t)\hat{\mathcal{X}}_t) \in \cdot | \hat{\mathcal{X}}_t] \qquad \text{a.s.} \end{aligned}$$

$\forall t, r \geq 0, t + q < r$. Taking the limit $r \uparrow \infty$, we see that also (3.47) holds. $\square$

3.3.2. *Identification of the embedded tree.* Our next step is to identify the embedded tree $I$ in Proposition 39 as a binary splitting particle system. For $t \geq 0$, define equivalence relation $\overset{t-}{\sim}$ and $\overset{t+}{\sim}$ on $I$ by

$$(3.57) \quad \begin{aligned} w \overset{t-}{\sim} v & \quad \text{if and only if } \pi_{[0,t)}(w) = \pi_{[0,t)}(v), \\ w \overset{t+}{\sim} v & \quad \text{if and only if } \pi_{[0,t+\varepsilon]}(w) = \pi_{[0,t+\varepsilon]}(v) \text{ for some } \varepsilon > 0, \end{aligned}$$

and let $I_{t-}$ and $I_{t+}$ denote the collections of $\overset{t-}{\sim}$ and $\overset{t+}{\sim}$ equivalence classes in $I$, respectively. Define counting measures $\hat{X}_{t-}$ and $\hat{X}_{t+}$ on $\mathcal{D}_E[0,t]$ by

$$\hat{X}_{t-} := \sum_{w \in I_{t-}} \delta_{\pi_{[0,t]}(w)}, \qquad t > 0,$$

(3.58)
$$\hat{X}_t := \sum_{w \in I_{t+}} \delta_{\pi_{[0,t]}(w)}, \qquad t \geq 0.$$

It is not hard to see that $\hat{X} = (\hat{X}_t)_{t \geq 0}$ has right-continuous sample paths with left limits given by $\hat{X}_{t-}$ and that

(3.59) $$I_t = \hat{X}_t \qquad \text{a.s. } \forall t \geq 0.$$

Note that the a.s. and the $\forall t \geq 0$ cannot be interchanged here, since $\hat{X}_t$ is not a simple point measure at those (random) times when $|I_t| < |I_{t+}|$, that is, when splitting occurs.

LEMMA 40 (Identification of the embedded tree). *The process $\hat{X}$ is the $(G, \alpha, 0)$-particle system started at time $0$ in $\mathrm{Pois}(\mu)$.*

PROOF. By (3.59) and (3.47),

(3.60) $$P[\hat{X}_t \in \cdot | (\hat{\mathcal{X}}_s)_{0 \leq s \leq t}] = P[\mathrm{Pois}(\hat{\mathcal{X}}_t) \in \cdot | \hat{\mathcal{X}}_t] \qquad \text{a.s. } \forall t \geq 0.$$

Let $\hat{X}'$ denote the $(G, \alpha, 0)$-particle system started at time $0$ in $\mathrm{Pois}(\mu)$. The time-inhomogeneous log-Laplace semigroup $(\hat{\mathcal{U}}_{s,t})_{0 \leq s \leq t}$ of the historical $(G, \alpha, \alpha)$-superprocess $\hat{\mathcal{X}}$ and the time-inhomogeneous generating semigroup $(\hat{U}_{s,t})_{0 \leq s \leq t}$ of the historical $(G, \alpha, 0)$-particle system $\hat{X}'$ are defined by the same Cauchy integral equation. Hence

(3.61) $$\hat{\mathcal{U}}_{s,t} f = \hat{U}_{s,t} f, \qquad 0 \leq s \leq t, f \in B_{[0,1]}(\mathcal{D}_E[0,t]).$$

Therefore, we may reason exactly as in the proof of Lemma 1 to see that

(3.62) $$P^{0, \mathrm{Pois}(\mu)}[\hat{X}'_t \in \cdot] = P^{0,\mu}[\mathrm{Pois}(\hat{\mathcal{X}}_t) \in \cdot], \qquad t \geq 0, \mu \in \mathcal{M}(E).$$

Combining (3.60) and (3.62), we see that

(3.63) $$P[\hat{X}_t \in \cdot] = P[\hat{X}'_t \in \cdot], \qquad t \geq 0.$$

It follows from our definition of $\hat{X}$ that

(3.64) $$\hat{X}_s = \mathrm{supp}\,(\hat{X}_t \circ \pi_{[0,s]}^{-1}) \qquad \text{a.s. } \forall 0 \leq s \leq t.$$

By a straightforward analogue of Lemma 23(a) for historical particle systems, $\mathrm{supp}(\hat{X}'_t \circ \pi_{[0,s]}^{-1}) \subset \mathrm{supp}(\hat{X}'_s)$ a.s. $\forall 0 \leq s \leq t$. Since the death rate of $\hat{X}'$ is zero, particles cannot become extinct and, therefore, in fact $\mathrm{supp}(\hat{X}'_t \circ \pi_{[0,s]}^{-1}) = \mathrm{supp}(\hat{X}'_s)$ a.s. $\forall 0 \leq s \leq t$. Since $\hat{X}'_s$ is a.s. a simple point measure [which follows from (3.63) and the fact that $\hat{X}_s$ is a.s. a simple point measure], $X'$ satisfies, in analogy with (3.64),

(3.65) $$\hat{X}'_s = \mathrm{supp}\,(\hat{X}'_t \circ \pi_{[0,s]}^{-1}) \qquad \text{a.s. } \forall 0 \leq s \leq t.$$

It follows from (3.63)–(3.65) that

(3.66) $P[(\hat{X}_{t_1},\ldots,\hat{X}_{t_n}) \in \cdot] = P[(\hat{X}'_{t_1},\ldots,\hat{X}'_{t_n}) \in \cdot], \qquad 0 \leq t_1 < t_2 < \cdots < t_n.$

Since $\hat{X}$ and $\hat{X}'$ have right-continuous sample paths, $\hat{X}$ and $\hat{X}'$ are equal in distribution. □

3.3.3. *Proof of the main theorems.* Theorems 6, 8 and 9 can be combined into the following theorem.

THEOREM 41 (Main results). *Let $\hat{\mathcal{X}}$ be the historical $(G,\alpha,\beta)$-superprocess started at time 0 in $\mu \in \mathcal{M}(E)$. Assume that $h \in \mathcal{D}(G)$ satisfies $h > 0$ and, for some $\gamma \in \mathcal{C}_+(E)$,*

(3.67) $$Gh + \beta h - \alpha h^2 = -\gamma h.$$

*Then $\hat{\mathcal{X}}$ can be coupled to the historical $(G^h, h\alpha, \gamma)$-particle system $\hat{X}$ started in $\hat{X}_0 = \mathrm{Pois}(h\mu)$ such that*

(3.68) $P[\hat{X}_t \in \cdot | (\hat{\mathcal{X}}_s)_{0 \leq s \leq t}] = P[\mathrm{Pois}((h \circ \pi_t)\hat{\mathcal{X}}_t) \in \cdot | \hat{\mathcal{X}}_t]$ *a.s.* $\forall t \geq 0.$

*If, in addition, $\mathcal{U} = \mathcal{U}(G,\alpha,\beta)$ satisfies $\sup_{x \in E} \mathcal{U}_t \infty(x) < \infty$ for some $t > 0$, then $p := \lim_{t \uparrow \infty} \mathcal{U}_t \infty \leq h$ and the coupling may be chosen such that, moreover,*

(3.69) $\mathrm{supp}(\hat{X}_t) \supset \mathrm{supp}(\hat{\mathcal{X}}_r \circ \pi_{[0,t]}^{-1}),$ *r-eventually $\forall t \geq 0$ a.s.*

*If, in addition, $\gamma = 0$, then $p = h$ and the coupling may be chosen such that equality holds r-eventually in* (3.69).

PROOF. Under the additional assumptions that (i) $\mu$ is atomless and the Feller process with generator $G$ has the distinct path property, (ii) $\gamma = 0$ and (iii) $h = 1$, the statement follows from Proposition 39 and Lemma 40. We now remove these assumptions one by one.

(i) *Generalization to measures with atoms.* Let $\eta$ be a Feller process in a compact metrizable space $F$ such that $\eta$ has the distinct path property (e.g., Brownian motion on the unit circle). Let $G'$ denote the generator of the Feller process $(\xi,\eta)$ in $E \times F$, where for given initial conditions, $\xi$ and $\eta$ evolve independently. Put $\alpha'(x,y) := \alpha(x)$ and $\beta'(x,y) := \beta(x)$. Let $\psi_t$ denote the projection from $\mathcal{D}_{E \times F}[0,t]$ to $\mathcal{D}_E[0,t]$. Let $\hat{\mu}$ and $\hat{\rho}$ be finite measures on $\mathcal{D}_E[0,s]$ and $\mathcal{D}_F[0,s]$, respectively, and assume that $\hat{\rho}$ is atomless. If $\hat{\mathcal{X}}'$ is the historical $(G',\alpha',\beta')$-superprocess started at time $s$ in $\hat{\mu} \otimes \hat{\rho}$, then, by Lemma 18,

(3.70) $$\hat{\mathcal{X}}_t := \hat{\mathcal{X}}'_t \circ \psi_t^{-1}, \qquad t \geq s$$

is the historical $(G, \alpha, \beta)$-superprocess started at time $s$ in $\hat{\mu}$. Moreover, $\hat{\mathcal{X}}'_t$ is atomless a.s. $\forall t \geq s$ and its underlying motion has the distinct path property. The statements for $\hat{\mathcal{X}}$ now follow from the statements for $\hat{\mathcal{X}}'$ by projection.

(ii) *Generalization to* $\gamma \neq 0$. Note that since we are still assuming $h = 1$, (3.67) reduces to $\alpha - \beta = \gamma$. Set $E^\dagger := E \cup \{\dagger\}$, where $\dagger$ is an isolated cemetery point that does not belong to $E$. Define a linear operator $G^\dagger$ on $\mathcal{C}(E^\dagger)$ by

(3.71)
$$G^\dagger f(x) := Gf(x) + \gamma(x)(f(\dagger) - f(x)), \qquad x \in E,$$
$$G^\dagger f(\dagger) := 0,$$

where $\mathcal{D}(G^\dagger)$ consists of those $f \in \mathcal{C}(E^\dagger)$ such that the restriction of $f$ to $E$ is in $\mathcal{D}(G)$. Set, moreover,

(3.72)
$$\alpha^\dagger(x) := \alpha(x), \qquad x \in E,$$
$$\alpha^\dagger(\dagger) := 1.$$

Let $\hat{\mathcal{X}}^\dagger$ denote the historical $(G^\dagger, \alpha^\dagger, \alpha^\dagger)$-superprocess started at time 0 in $\mu \in \mathcal{M}(E)$ and let $\hat{X}^\dagger$ denote the historical $(G^\dagger, \alpha^\dagger, 0)$-particle system started at time 0 in $\text{Pois}(\mu)$. For $t \geq 0$, let $\hat{\mathcal{X}}_t$ and $\hat{X}_t$ denote the restrictions of $\hat{\mathcal{X}}_t^\dagger$ and $\hat{X}_t^\dagger$ to $\mathcal{D}_E[0, t]$, respectively. Elementary considerations involving the log-Laplace semigroups of $\hat{\mathcal{X}}_t^\dagger$ and $\hat{X}_t^\dagger$ show that $(\hat{\mathcal{X}}_t)_{t \geq 0}$, so defined, is the historical $(G, \alpha, \beta)$-superprocess, and that $(\hat{X}_t)_{t \geq 0}$ is the historical $(G, \alpha, \gamma)$-particle system.

By what we have already proved, $\hat{\mathcal{X}}^\dagger$ and $\hat{X}^\dagger$ may be coupled such that

(3.73) $\quad P[\hat{X}_t^\dagger \in \cdot | (\hat{\mathcal{X}}_s^\dagger)_{0 \leq s \leq t}] = P[\text{Pois}(\hat{\mathcal{X}}_t^\dagger) \in \cdot | \hat{\mathcal{X}}_t^\dagger] \qquad \text{a.s. } \forall t \geq 0,$

which implies (3.68). If, in addition, $\sup_{x \in E} \mathcal{U}_t \infty(x) < \infty$ for some $t > 0$, then using the fact that $\alpha^\dagger(\dagger) = 1$, it is not hard to show that also $\sup_{x \in E^\dagger} \mathcal{U}_t^\dagger \infty(x) < \infty$ for some $t > 0$ and, therefore, by what we have already proved,

$$p^\dagger := \lim_{t \uparrow \infty} \mathcal{U}_t^\dagger \infty = 1$$

and the coupling between $\hat{\mathcal{X}}^\dagger$ and $\hat{X}^\dagger$ may be chosen such that, moreover,

(3.74) $\quad \text{supp}(\hat{X}_t^\dagger) = \text{supp}(\hat{\mathcal{X}}_r^\dagger \circ \pi_{[0,t]}^{-1}), \qquad r\text{-eventually } \forall t \geq 0 \text{ a.s.}$

By Lemma 19(b) and the fact that $\dagger$ is a trap for the underlying motion, $\hat{\mathcal{X}}^\dagger$ is concentrated on paths that are trapped in $\dagger$, once they reach $\dagger$ and, therefore,

(3.75)
$$\text{supp}(\hat{X}_t) = \text{supp}(\hat{X}_t^\dagger) \cap \mathcal{D}_E[0, t]$$
$$= \text{supp}(\hat{\mathcal{X}}_r^\dagger \circ \pi_{[0,t]}^{-1}) \cap \mathcal{D}_E[0, t] \supset \text{supp}(\hat{\mathcal{X}}_r \circ \pi_{[0,t]}^{-1})$$

$\forall\, 0 \leq t \leq r$ a.s. Formulas (3.74) and (3.75) imply (3.69). Finally, for all $x \in E$,

$$
\begin{aligned}
p(x) &= -\log P^{\delta_x}[\mathcal{X}_t = 0 \text{ eventually}] \\
&\leq -\log P^{\delta_x}[\mathcal{X}_t^\dagger = 0 \text{ eventually}] = p^\dagger(x) = 1.
\end{aligned}
\tag{3.76}
$$

(iii) *Generalization to* $h \neq 1$. Set $\hat{\mathcal{X}}_t^h(dw) := h(w_t)\hat{\mathcal{X}}_t(dw)$, $t \geq 0$. By Lemma 5, $\hat{\mathcal{X}}^h$ is the historical $(G^h, \alpha^h, \beta^h)$-superprocess, where $G^h$ is defined in (1.16) and $\alpha^h := h\alpha$, $\beta^h := \beta + \frac{Gh}{h}$. Formula (3.67) implies that

$$-\gamma = \beta^h - \alpha^h \leq 0. \tag{3.77}$$

Therefore the statements follow from what we have already proved. $\square$

**Acknowledgments.** We are grateful to Tom Kurtz for providing a simpler proof of Lemma 3 and for referring us to some references such as [15]. We also thank an anonymous referee for a very careful reading, which led to a better exposition, and for suggesting the improved proof of Lemma 31. We thank János Engländer and Anita Winter for helpful comments.

WEIERSTRASS INSTITUTE FOR
APPLIED ANALYSIS AND STOCHASTICS
MOHRENSTRASSE 39
D-10117 BERLIN
GERMANY
E-MAIL: fleischm@wias-berlin.de
URL: www.wias-berlin.de/fleischm

UNIVERSITY ERLANGEN–NUREMBERG
BISMARCKSTRASSE 1 1/2
91054 ERLANGEN
GERMANY
E-MAIL: swart@mi.uni-erlangen.de
URL: www.mi.uni-erlangen.de/swart